# BU-CE-1201

# Technical Report on Hypergraph-Partitioning-Based Models and Methods for Exploiting Cache Locality in Sparse-Matrix Vector Multiplication

Kadir Akbudak, Enver Kayaaslan and Cevdet Aykanat

February 2012

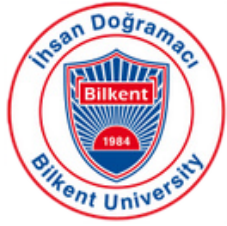

# TECHNICAL REPORT ON HYPERGRAPH-PARTITIONING-BASED MODELS AND METHODS FOR EXPLOITING CACHE LOCALITY IN SPARSE-MATRIX VECTOR MULTIPLICATION

KADIR AKBUDAK∗, ENVER KAYAASLAN ∗, AND CEVDET AYKANAT ∗

**Abstract.** The sparse matrix-vector multiplication (SpMxV) is a kernel operation widely used in iterative linear solvers. The same sparse matrix is multiplied by a dense vector repeatedly in these solvers. Matrices with irregular sparsity patterns make it difficult to utilize cache locality effectively in SpMxV computations. In this work, we investigate single- and multiple-SpMxV frameworks for exploiting cache locality in SpMxV computations. For the single-SpMxV framework, we propose two cache-size-aware top-down row/column-reordering methods based on 1D and 2D sparse matrix partitioning by utilizing the column-net and enhancing the row-column-net hypergraph models of sparse matrices. The multiple-SpMxV framework depends on splitting a given matrix into a sum of multiple nonzero-disjoint matrices so that the SpMxV operation is performed as a sequence of multiple input- and output-dependent SpMxV operations. For an effective matrix splitting required in this framework, we propose a cache-size-aware top-down approach based on 2D sparse matrix partitioning by utilizing the row-column-net hypergraph model. For this framework, we also propose two methods for effective ordering of individual SpMxV operations. The primary objective in all of the three methods is to maximize the exploitation of temporal locality. We evaluate the validity of our models and methods on a wide range of sparse matrices using both cache-miss simulations and actual runs by using OSKI. Experimental results show that proposed methods and models outperform state-of-the-art schemes.

**1. Introduction.** Sparse matrix-vector multiplication (SpMxV) is an important kernel operation in iterative linear solvers used for the solution of large, sparse, linear systems of equations. In these iterative solvers, the SpMxV operation $y \leftarrow Ax$ is repeatedly performed with the same large, irregularly sparse matrix $A$. Irregular access pattern during these repeated SpMxV operations causes poor usage of CPU caches in today's deep memory hierarchy technology. However, SpMxV operation has a potential to exhibit very high performance gains if temporal and spatial localities are respected and exploited properly. Here, temporal locality refers to the reuse of data words (e.g., $x$-$vector$ entries) within relatively small time durations, whereas spatial locality refers to the use of data words (e.g., matrix nonzeros) within relatively close storage locations (e.g., in the same lines).

In this work, we investigate two distinct frameworks for the SpMxV operation: single-SpMxV and multiple-SpMxV frameworks. In the single-SpMxV framework, the $y$-$vector$ results are computed by performing a single SpMxV operation $y \leftarrow Ax$. In the multiple-SpMxV framework, $y \leftarrow Ax$ operation is computed as a sequence of multiple input- and output-dependent SpMxV operations, $y \leftarrow y + A^k x$ for $k = 1, \ldots, k$, where $A = A^1 + \cdots + A^K$. For the single-SpMxV framework, we propose two cache-size-aware row/column reordering methods based on top-down 1D and 2D partitioning of a given sparse matrix. The 1D-partitioning-based method relies on transforming a sparse matrix into a singly-bordered block-diagonal (SB) form by utilizing the column-net hypergraph model [7–9]. The 2D-partitioning-based method relies on transforming a sparse matrix into a doubly-bordered block-diagonal (DB) form by utilizing the row-column-net hypergraph model [7, 11]. We provide upper bounds on the number of cache misses based on these transformations, and show that the objectives in the transformations based on partitioning the respective hypergraph models correspond to minimizing these upper bounds. In the 1D-partitioning-based method, the column-net hypergraph model correctly encapsulates the minimization of the respective upper bound. For the 2D-partitioning-based method, we propose an enhancement to the row-column-net hypergraph model to encapsulate the minimization of the respective upper bound on the number of cache misses. The primary objective in both methods is to maximize the exploitation of the temporal locality due to the access of $x$-$vector$ entries, whereas exploitation of the spatial locality due to the access of $x$-$vector$ entries is a secondary objective. In this paper, we claim that exploiting temporal locality is more important than exploiting spatial locality (for practical line sizes) in SpMxV operations that involve irregularly sparse matrices.

---

∗Computer Engineering Department, Bilkent University, Ankara, Turkey (kadir@cs.bilkent.edu.tr, enver@cs.bilkent.edu.tr, aykanat@cs.bilkent.edu.tr)

The multiple-SpMxV framework depends on splitting a given matrix into a sum of multiple nonzero-disjoint matrices so that the SpMxV operation is computed as a sequence of multiple SpMxV operations. For an effective matrix splitting required in this framework, we propose a cache-size-aware top-down approach based on 2D sparse matrix partitioning by utilizing the row-column-net hypergraph model [7, 11]. We provide an upper bound on the number of cache misses based on this matrix-splitting, and show that the objective in the hypergraph-partitioning (HP) based matrix partitioning exactly corresponds to minimizing this upper bound. The primary objective in this method is to maximize the exploitation of the temporal locality due to the access of both $x$-$vector$ and $y$-$vector$ entries. For this framework, we also propose two methods for effective ordering of individual SpMxV operations.

We evaluate the validity of our models and methods on a wide range of sparse matrices. The experiments are carried out in two different settings: cache-miss simulations and actual runs by using OSKI (BeBOP Optimized Sparse Kernel Interface Library) [41]. Experimental results show that the proposed methods and models outperform state-of-the-art schemes and also these results conform to our expectation that temporal locality is more important than spatial locality in SpMxV operations that involve irregularly sparse matrices.

The rest of the paper is organized as follows: Background material is introduced in Section 2. In Section 3, we review some of the previous works about iteration/data reordering and matrix transformations for exploiting locality. The two frameworks along with our contributed models and methods are described in Sections 4 and 5. We present the experimental results in Section 6. Finally, the paper is concluded in Section 7.

**2. Background.** Several sparse-matrix storage schemes utilized in SpMxV are summarized in Section 2.1. Data locality issues during the SpMxV operation are discussed in Section 2.2. Section 2.3 summarizes the HP problem, whereas Section 2.4 discusses hypergraph models and methods for sparse-matrix partitioning. Finally, bipartite graph model for sparse matrices is given in Section 2.5.

**2.1. Sparse-matrix storage schemes.** There are two standard sparse-matrix storage schemes for the SpMxV operation: *Compressed Storage by Rows* (CSR) and *Compressed Storage by Columns* (CSC) [14, 33]. In this paper, we restrict our focus on the conventional SpMxV operation using the CSR storage scheme without loss of generality, whereas cache aware techniques such as prefetching, blocking, etc. are out of the scope of this paper. In the following paragraphs, we review the standard CSR scheme and two CSR variants.

The CSR scheme contains three 1D arrays: *nonzero*, *colIndex* and *rowStart*. The values and the column indices of nonzeros are respectively stored in row-major order in the $nonzero$ and $colIndex$ arrays in a one-to-one manner. That is, $colIndex[k]$ stores the column index of the nonzero stored in $nonzero[k]$. The *rowStart* array stores the index of the first nonzero of each row in the *nonzero* and $colIndex$ arrays. Algorithm 1 shows SpMxV utilizing the CSR storage scheme for an $m \times n$ sparse matrix. Each outer for-loop iteration of Algorithm 1 corresponds to the inner product of the respective sparse row with the dense input vector $x$.

The *zig-zag CSR* (ZZCSR) scheme is proposed to reduce end-of-row cache misses [43]. In ZZCSR, nonzeros are stored in increasing column-index order in even-numbered rows, whereas they are stored in decreasing index order in odd-numbered rows, or vice versa.

**Algorithm 1** SpMxV using CSR scheme

**Require:** $nonzero$, $colIndex$ and $rowStart$ arrays of an $m \times n$ sparse matrix $A$
   a dense input vector $x$
**Output:** dense vector $y$
 1: **for** $i \leftarrow 1$ *to* $m$ **do**
 2:   $sum \leftarrow 0.0$
 3:   **for** $k \leftarrow rowStart\,[i]$ **to** $rowStart\,[i+1] - 1$ **do**
 4:     $sum \leftarrow sum + nonzero[k] * x[\,colIndex\,[k]]$
 5:   **end for**
 6:   $y[i] \leftarrow sum$
 7: **end for**
 8: **return** y

**Algorithm 2** SpMxV using ICSR scheme [27]

**Require:** $nonzero$, $colDiff$ and $rowJump$ arrays of an $m \times n$ sparse matrix $A$ with $nnz$ nonzeros,
   a dense input vector $x$
**Output:** dense vector $y$
 1: $i \leftarrow rowJump$
 2: $j \leftarrow colDiff\,[0]$
 3: $k \leftarrow 0$
 4: $r \leftarrow 1$
 5: $sum \leftarrow 0.0$
 6: **for** $k \leftarrow 1$ **to** $nnz$ **do**
 7:   $sum \leftarrow sum + nonzero[k] * x[j]$
 8:   $k \leftarrow k + 1$
 9:   $j \leftarrow j + colDiff\,[k]$
10:   **if** $j \geq n$ **then**
11:     $y[i] \leftarrow sum$
12:     $sum \leftarrow 0.0$
13:     $j \leftarrow j - n$
14:     $i \leftarrow i + rowJump\,[r]$
15:     $r \leftarrow r + 1$
16:   **end if**
17: **end for**
18: **return** y

The *Incremental Compressed Storage by Rows* (ICSR) scheme [27] which is given in Algorithm 2s reported to decrease instruction overhead by using pointer arithmetic. In ICSR, the $colIndex$ array is replaced with the $colDiff$ array, which stores the increments in the column indices of the successive nonzeros stored in the nonzero array. The $rowStart$ array is replaced with the $rowJump$ array which stores the increments in the row indices of the successive nonzero rows. The beginning of a new row is signalled by causing an increment value $j$ to overflow $n$ so that $j - n$ shows the column index of the first nonzero in the next row. For this purpose, nonzeros of each row are stored in increasing column-index order. The ICSR scheme has the advantage of handling zero rows efficiently since it avoids the use of the $rowStart$ array. This feature of ICSR is exploited in our multiple-SpMxV framework since this scheme introduces many zero rows in the individual sparse matrices.

**2.2. Data locality in CSR-based SpMxV.** In accessing matrix nonzeros, temporal locality is not feasible since the elements of each of the $nonzero$, $colIndex$ ($colDiff$ in ICSR) and $rowStart$ ($rowJump$ in ICSR) arrays are accessed only once. Spatial locality is feasible and it is achieved automatically by nature of the CSR scheme since the elements of each of the three arrays are accessed consecutively.

In accessing $y$-$vector$ entries, temporal locality is not feasible since each $y$-$vector$ result is written only once to the memory. As a different view, temporal locality can be considered as feasible but automatically achieved especially at the register level because of the summation of scalar nonzero and $x$-$vector$ entry product results to the temporary variable $sum$. Spatial locality is feasible and it is achieved automatically since the $y$-$vector$ entry results are stored consecutively.

In accessing $x$-$vector$ entries, both temporal and spatial localities are feasible. Temporal locality is feasible since each $x$-$vector$ entry may be accessed multiple times. However, exploiting the temporal and spatial localities for the $x$-$vector$ is the major concern in the CSR scheme since $x$-$vector$ entries are accessed through a $colIndex$ array ($colDiff$ in ICSR) in a non-contiguous and irregular manner.

These locality issues can be solved by reordering rows/columns of matrix $A$ and the exploitation level of these data localities depends both on the existing sparsity pattern of matrix $A$ and the effectiveness of reordering heuristics.

**2.3. Hypergraph partitioning.** A hypergraph $\mathcal{H} = (\mathcal{V}, \mathcal{N})$ is defined as a set $\mathcal{V}$ of vertices and a set $\mathcal{N}$ of nets (hyperedges). Every net $n_j \in \mathcal{N}$ connects a subset of vertices, i.e., $n_j \subseteq \mathcal{V}$. Weights and costs can be associated with vertices and nets, respectively. We use $w(v_i)$ to denote the weight of vertex $v_i$ and $cost(n_j)$ to denote the cost of net $n_j$. Given a hypergraph $\mathcal{H} = (\mathcal{V}, \mathcal{N})$, $\{\mathcal{V}_1, \ldots, \mathcal{V}_K\}$ is called a $K$-way partition of the vertex set $\mathcal{V}$ if vertex parts are mutually disjoint and exhaustive. A $K$-way vertex partition of $\mathcal{H}$ is said to satisfy the partitioning constraint if $W_k \leq W_{avg}(1+\varepsilon)$ for $k = 1, 2, \ldots, K$. $W_k$ denotes weight of a part $\mathcal{V}_k$ and is defined as the sum of weights of vertices in $\mathcal{V}_k$. $W_{avg}$ is the average part weight and $\varepsilon$ represents a predetermined, maximum allowable imbalance ratio.

In a partition of $\mathcal{H}$, a net that connects at least one vertex in a part is said to *connect* that part. *Connectivity* $\lambda(n_j)$ of a net $n_j$ denotes the number of parts connected by $n_j$. A net $n_j$ is said to be *cut* if it connects more than one part (i.e., $\lambda(n_j) > 1$), and *uncut* (*internal*) otherwise (i.e., $\lambda(n_j) = 1$). The set of cut nets of a partition is denoted as $\mathcal{N}_{cut}$. The partitioning objective is to minimize the cutsize defined over the cut nets. There are various cutsize definitions. Two relevant cutsize definitions are the cut-net and connectivity [6] metrics.

$$cutsize-cutnet = \sum_{n_j \in \mathcal{N}_{cut}} cost(n_j) \qquad cutsize-con = \sum_{n_j \in \mathcal{N}_{cut}} \lambda(n_j)\, cost(n_j) \qquad (2.1)$$

In the cut-net metric, each cut net $n_j$ incurs $cost(n_j)$ to the cutsize, whereas in the connectivity metric, each cut net incurs $\lambda(n_j)\, cost(n_j)$ to the cutsize. The HP problem is known to be NP-hard [28]. There exists several successful HP tools such as hMeTiS [26], PaToH [10] and Mondriaan [40], all of which apply the multilevel framework.

The *recursive bisection* (RB) paradigm is widely used in $K$-way HP and known to be amenable to produce good solution qualities. In the RB paradigm, first, a 2-way partition of the hypergraph is obtained. Then, each part of the bipartition is further bipartitioned in a recursive manner until the desired number $K$ of parts is obtained or part weights drop below a given part-weight threshold $W_{max}$. In RB-based HP, the cut-net removal and cut-net splitting schemes [9] are used to capture the cut-net and connectivity cutsize metrics, respectively. The RB paradigm is inherently suitable for partitioning hypergraphs when K is not known in advance. Hence, the RB paradigm can be successfully utilized in clustering rows/columns for cache-size-aware row/column reordering.

**2.4. Hypergraph models for sparse matrix partitioning.** Recently, several successful hypergraph models are proposed for partitioning sparse matrices [7, 9]. The relevant ones are row-net, column-net and row-column-net (fine-grain) models.

In the row-net hypergraph model [7–9] $\mathcal{H}_{\text{RN}}(A) = (\mathcal{V}_\mathcal{C}, \mathcal{N}_\mathcal{R})$ of matrix $A$, there exist one vertex $v_j \in \mathcal{V}_\mathcal{C}$ and one net $n_i \in \mathcal{N}_\mathcal{R}$ for each column $c_j$ and row $r_i$, respectively. The weight $w(v_j)$ of a vertex $v_j$ is set to the number of nonzeros in column $c_j$. The net $n_i$ connects the vertices corresponding to the columns that have a nonzero entry in row $r_i$. Every net $n_i \in \mathcal{N}_\mathcal{R}$ has unit cost, i.e., $cost(n_i) = 1$. In the column-net hypergraph model [7–9] $\mathcal{H}_{CN}(A) = (\mathcal{V}_\mathcal{R}, \mathcal{N}_\mathcal{C})$ of matrix $A$, there exist one vertex $v_i \in \mathcal{V}_\mathcal{R}$ and one net $n_j \in \mathcal{N}_\mathcal{C}$ for each row $r_i$ and column $c_j$, respectively. The weight $w(v_i)$ of a vertex $v_i$ is set to the number of nonzeros in row $r_i$. Net $n_j$ connects the vertices corresponding to the rows that have a nonzero entry in column $c_j$. Every net $n_j$ has unit cost, i.e., $cost(n_j) = 1$.

In the row-column-net model [11] $\mathcal{H}_{RCN}(A) = (\mathcal{V}_\mathcal{Z}, \mathcal{N}_{\mathcal{RC}})$ of matrix $A$, there exists one vertex $v_{ij} \in \mathcal{V}_\mathcal{Z}$ corresponding to each nonzero $a_{ij}$ in matrix $A$. In net set $\mathcal{N}_{\mathcal{RC}}$, there exists a row-net $n_i^r$ for each row $r_i$, and there exists a column-net $n_j^c$ for each column $c_j$. Every row net and column net have unit cost. Row-net $n_i^r$ connects the vertices corresponding to the nonzeros in row $r_i$, and column-net $n_j^c$ connects the vertices corresponding to the nonzeros in column $c_j$. Note that each vertex is connected by exactly two nets.

In [3], it is shown that row-net and column-net models can also be used for transforming a sparse matrix into a $K$-way SB form through row and column reordering. In particular, the row-net model can be used for permuting a matrix into a rowwise SB form, whereas the column-net model can be used for permuting a matrix into a columnwise SB form. Here we will briefly describe how a $K$-way partition of the column-net model can be decoded as a row/column reordering for this purpose and a dual discussion holds for the row-net model.

A $K$-way vertex partition $\{\mathcal{V}_1, \ldots, \mathcal{V}_K\}$ of $\mathcal{H}_{CN}(A)$ is considered as inducing a $(K+1)$-way partition $\{\mathcal{N}_1, \ldots, \mathcal{N}_K; \mathcal{N}_{cut}\}$ on the net set of $\mathcal{H}_{CN}(A)$. Here $\mathcal{N}_k$ denotes the set of internal nets of vertex part $\mathcal{V}_k$, whereas $\mathcal{N}_{cut}$ denotes the set of cut nets. The vertex partition is decoded as a partial row reordering of matrix $A$ such that the rows associated with vertices in $\mathcal{V}_{k+1}$ are ordered after the rows associated with vertices $\mathcal{V}_k$, for $k = 1, \ldots, K-1$. The net partition is decoded as a partial column reordering of matrix $A$ such that the columns associated with nets in $\mathcal{N}_{k+1}$ are ordered after the columns associated with nets in $\mathcal{N}_k$, for $k = 1, \ldots, K-1$, whereas the columns associated with the cut nets are ordered last to constitute the column border.

**2.5. Bipartite graph model for sparse matrices.** In the *bipartite graph model* $\mathcal{B}(A) = (\mathcal{V}, \mathcal{E})$ of matrix $A$, there exists one row vertex $v_i^r \in \mathcal{R}$ representing row $r_i$, and there exists one column vertex $v_j^c \in \mathcal{C}$ representing column $c_j$, where $\mathcal{R}$ is the set of row vertices and $\mathcal{C}$ is the set of column vertices. These vertex sets $\mathcal{R}$ and $\mathcal{C}$ form the vertex bipartition $\mathcal{V} = \mathcal{R} \cup \mathcal{C}$. There is an edge between vertices $v_i^r \in \mathcal{R}$ and $v_j^c \in \mathcal{C}$ if and only if the respective matrix entry $a_{ij}$ is nonzero.

**3. Related work.** The main focus of this work is to perform iteration and data reordering, without changing the conventional CSR-based SpMxV codes, whereas cache aware techniques such as prefetching, blocking, etc. are out of the scope of this paper. So we summarize the related work on iteration and data reordering for irregular applications which usually use index arrays to access other arrays. Iteration and data reordering approaches can also be categorized as dynamic and static. Dynamic schemes [12, 13, 15, 19, 35] achieve runtime reordering transformations by analyzing the irregular memory access patterns through adopting inspector/executor strategy [29]. Reordering rows/columns of irregularly sparse matrices to exploit locality during SpMxV operations can be considered as a static case of such general iteration/data reordering problem. We call it a static case [32, 38, 42, 43] since the sparsity pattern of matrix $A$ together with the CSR- or CSC-based SpMxV scheme determines the memory access pattern. In the CSR scheme, iteration order corresponds to row order of matrix $A$ and data order corresponds to column order, whereas a dual discussion applies for

CSC.

Dynamic and static transformation heuristics mainly differ in the preprocessing times. Fast heuristics are usually used for dynamic reordering transformations, whereas much more sophisticated heuristics are used for static case. The preprocessing time for the static case can amortize the performance improvement during repeated computations with the same memory access pattern. Repeated SpMxV computations involving the same matrix or matrices with the same sparsity pattern constitute a very typical case of such static case.

Ding and Kennedy [15] propose the locality grouping and consecutive packing (CPACK) heuristics for runtime iteration and data reordering, respectively. The locality grouping heuristic traverses the data objects in a given order and clusters all the iterations that access the first data item, then the second, and etc. The CPACK heuristic reorders the data objects on a first-touch-first basis. The locality grouping heuristic is also referred to as consecutive packing for iterations (CPACKIter) in [35] and this heuristic is equivalent to the iteration reordering heuristic proposed by Das et al. [13] As also mentioned in [15, 19], these heuristics suffer from not explicitly considering different reuse patterns of different data objects because the data objects and iterations are traversed in a given order.

Space-filling curves such as Hilbert and Morton as well as recursive storage schemes such as quadtree are used for iteration reordering in improving locality in dense matrix operations [16, 17, 25] and in sparse matrix operations [18]. Space-filling curves [12] and hierarchical graph clustering algorithms (GPART) [19] are utilized for data reordering in improving locality in n-body simulation applications.

Strout et al. [34] integrate run-time data and iteration reordering transformations such as lexicographically grouping, CPACK and GPART into a compile time framework and they show that sparse tiling may improve performance of these transformations depending on the underlying architecture. Strout and Hovland [35] extend the work in [34] and propose hypergraph-based models for data and iteration reordering transformations. They introduce a temporal locality hypergraph model for ordering iterations to exploit temporal locality. They also generalize spatial locality graph model to spatial locality hypergraph model to encompass the applications having multiple arrays that are accessed irregularly. Additionally, they propose a modified algorithm like Breadth-First Search (BFS) for ordering data and iterations simultaneously, whereas Breadth-First Search is used for only data ordering in [2]. Strout and Hovland [35] also propose metrics to determine which reordering heuristic is expected to yield better performance.

Das et al. [13] use reordering techniques in their implementation of three-dimensional unstructured grid Euler-solver to improve cache utilization. They reorder unstructured mesh edges incident on the same node consecutively. They also use Reverse Cuthill McKee (RCM) method to reorder nodes of the mesh. Burgess and Giles [5] examine effects of reordering techniques in unstructured grid applications. They report that reordering meshes that are generated without any cache optimization may result increase in performance according to application: Original orderings give better results in Jacobi solver, whereas reordered meshes give better results in conjugate gradients method.

Al-Furaih and Ranka [2] introduce interaction graph model to investigate optimizations for unstructured iterative applications in which the computational structure remains static or changes only slightly through iterations. They compare several methods to reorder data elements through reordering the vertices of the interaction graph. They report that BFS, as a fast reordering heuristic, can be applied to a static structure once or to a dynamic structure between tens of iterations. The other reordering methods are based on top-down graph partitioning, BFS ordering after graph partitioning and reordering via finding connected components that can fit into the cache.

In the rest of this section, we discuss the related work on improving locality in SpMxV operations. Agarwal et al. [1] try to improve SpMxV by extracting dense block structures. Their methods consist of examining row blocks to find dense subcolumns and reorder these subcolumns consecutively. Temam and Jalby [36] analyze the cache miss behaviour of SpMxV. They report that cache hit ratio decreases as bandwidth of sparse matrix increases beyond the cache size and conclude that bandwidth reduction algorithms improve cache utilization.

Toledo [38] compares several techniques to reduce cache misses in SpMxV. He uses graph theoretic methods such as Cuthill McKee (CM), RCM and top-down graph partitioning for reordering matrices and other improvement techniques such as blocking, prefetching and instruction-level-related optimization. He reports that SpMxV performance cannot be improved through row/column reordering. White and Sadayappan [42] discuss data locality issues in SpMxV in detail. They compare SpMxV performance of CSR, CSC and blocked versions of CSR and CSC. They also propose a graph-partitioning-based row/column reordering method which is similar to that of Toledo. They report that they can not achieve performance improvement over the original ordering as also reported by Toledo [38]. Haque and Hossain [20] propose a column reordering method based on Gray Code.

There are several works on row/column reordering based on similar TSP formulations. Heras et al. [23] define four distance functions for edge weighting depending on the similarity of sparsity patterns between row/columns. Pichel et al. [31] use TSP-based reordering and blocking technique to show improvements in both single processor performance and multicomputer performance. Pichel et al. [30] compare the performance of a number of reordering techniques which utilize TSP, top-down graph partitioning, RCM, Approximate Minimum Degree on simultaneous multithreading architectures. Pınar and Heath [32] propose a TSP-based column reordering for permuting nonzeros of a given matrix into contiguous blocks with the objective of decreasing the number of indirections in the CSR-based SpMxV. They compare the performance of their method to that of the RCM technique.

In a very recent work, Yzelman and Bisseling [43] propose a row/column reordering scheme based on partitioning row-net hypergraph representation of a given sparse matrix for CSR-based SpMxV. They achieve spatial locality on $x$-$vector$ entries by clustering the columns with similar sparsity pattern. They also exploit temporal locality for $x$-$vector$ entries by using zig-zag property of ZZCSR and ZZICSR schemes mentioned in Section 2.1.

**4. Single-SpMxV framework.** In this framework, the $y$-$vector$ results are computed by performing a single SpMxV operation, i.e., $y \leftarrow Ax$. The objective in this scheme is to reorder the columns and rows of matrix $A$ for maximizing the exploitation of temporal and spatial locality in accessing $x$-$vector$ entries. That is, the objective is to find row and column permutation matrices $P_r$ and $P_c$ so that $y \leftarrow Ax$ is computed as $\hat{y} \leftarrow \hat{A}\hat{x}$, where $\hat{A} = P_r A P_c$, $\hat{x} = x P_c$ and $\hat{y} = P_r y$. For the sake of simplicity of presentation, reordered input and output vectors $\hat{x}$ and $\hat{y}$ will be referred to as $x$ and $y$ in the rest of the paper.

Recall that temporal locality in accessing $y$-$vector$ entries is not feasible, whereas spatial locality is achieved automatically because $y$-$vector$ results are stored and processed consecutively. Reordering the rows with similar sparsity pattern nearby increases the possibility of exploiting temporal locality in accessing $x$-$vector$ entries. Reordering the columns with similar sparsity pattern nearby increases the possibility of exploiting spatial locality in accessing $x$-$vector$ entries. This row/column reordering problem can also be considered as a row/column clustering problem and this clustering process can be achieved in two distinct ways: top-down and bottom-up. In this section, we propose and discuss cache-size-aware top-down approaches based on 1D and 2D partitioning of a given matrix. Although a bottom-up approach based on hierarchical clustering of rows/columns with similar patterns is feasible, such a scheme is not discussed in this work.

In Sections 4.1 and 4.2, we present two theorems that give the guidelines for a "good" cache-size-aware row/column reordering based on 1D and 2D matrix partitioning. These theorems provide upper bounds on the number of cache misses due to the access of $x$-$vector$ entries in the SpMxV operation performed on sparse matrices in two special forms, namely SB and DB forms. In these theorems, $\Phi_x(A)$ denotes the number of cache misses due to the access of $x$-$vector$ entries in a CSR-based SpMxV operation to be performed on matrix $A$.

In the theorems given in Sections 4 and 5, fully associative cache is assumed, since misses in a fully associative cache are capacity misses and are not conflict misses. That is, each data line in the main memory can be placed to any empty line in the fully associative cache without causing a conflict miss. In these theorems, a matrix/submatrix is said to fit into the cache if the size of the CSR storage of the matrix/submatrix together with the associated $x$ and $y$ vectors/subvectors is smaller than the size of the cache.

**4.1. Row/column reordering based on 1D matrix partitioning.** We consider a row/column reordering which permutes a given matrix $A$ into a $K$-way columnwise SB form

$$\hat{A} = A_{SB} = P_r A P_c = \begin{bmatrix} A_{11} & & & & A_{1B} \\ & A_{22} & & & A_{2B} \\ & & \ddots & & \vdots \\ & & & A_{KK} & A_{KB} \end{bmatrix} = \begin{bmatrix} R_1 \\ R_2 \\ \vdots \\ R_K \end{bmatrix}$$

$$= \begin{bmatrix} C_1 & C_2 & \ldots & C_K & C_B \end{bmatrix}. \quad (4.1)$$

Here, $A_{kk}$ denotes the $k$th diagonal block of $A_{SB}$. $R_k = [0 \ldots 0 \ A_{kk} \ 0 \ldots 0 \ A_{kB}]$ denotes the $k$th row slice of $A_{SB}$, for $k = 1, \ldots, K$. $C_k = \begin{bmatrix} 0 \ldots 0 \ A_{kk}^T \ 0 \ldots 0 \end{bmatrix}^T$ denotes the $k$th column slice of $A_{SB}$, for $k = 1, \ldots, K$, and $C_B$ denotes the column border as follows

$$C_B = \begin{bmatrix} A_{1B} \\ A_{2B} \\ \vdots \\ A_{KB} \end{bmatrix}. \quad (4.2)$$

Each column in the border $C_B$ is called a *row-coupling column* or simply a *coupling column*. Let $\lambda(c_j)$ denote the number of submatrices $R_k$ that contain at least one nonzero of column $c_j$ of matrix $A_{SB}$, i.e.,

$$\lambda(c_j) = |\{R_k \in A_{SB} : c_j \in R_k\}|. \quad (4.3)$$

In other words, $\lambda(c_j)$ denotes the row-slice connectivity or simply connectivity of column $c_j$ in $A_{SB}$. In this notation, a column $c_j$ is a coupling column if $\lambda(c_j) > 1$. Here and hereafter, a submatrix notation is interchangeably used to denote both a submatrix and the set of non-empty rows/columns that belong to that matrix. For example, in (4.3), $R_k$ denotes both the $k$th row slice of $A_{SB}$ and the set of columns that belong to submatrix $R_k$.

The individual $y \leftarrow Ax$ can be equivalently represented as $K$ output-independent but input-dependent SpMxV operations, i.e., $y_k \leftarrow R_k x$ for $k = 1, \ldots, K$, where each submatrix $R_k$ is assumed to be stored in CSR scheme. These SpMxV operations are input dependent because of the $x$-$vector$ entries corresponding to the coupling columns.

THEOREM 4.1. *Given a $K$-way SB form $A_{SB}$ of matrix $A$ such that every submatrix $R_k$ fits into the cache, then we have*

$$\Phi_x(A_{SB}) \leq \sum_{c_j \in A_{SB}} \lambda(c_j) \tag{4.4}$$

*Proof.* Since each submatrix $R_k$ fits into the cache, for each $c_j \in R_k$, $x_j$ will be loaded to the cache at most once during the $y_k \leftarrow R_k x$ multiply. Therefore for a column $c_j$, the maximum number of cache misses that can occur due to the access of $x_j$ is bounded above by $\lambda(c_j)$. Note that this worst case happens when no cache-reuse occurs in accessing $x$-*vector* entries during successive $y_k \leftarrow R_k x$ operations implicitly performed in $y \leftarrow Ax$. □

Theorem 4.1 leads us to a cache-size-aware top-down row/column reordering through an $A$-to-$A_{SB}$ transformation that minimizes the upper bound given in (4.4) for $\Phi_x(A_{SB})$. Minimizing this sum relates to minimizing the number of cache misses due to the loss of temporal locality.

As discussed in [3], this $A$-to-$A_{SB}$ transformation problem can be formulated as an HP problem using the column-net model of matrix $A$ with the part size constraint of cache size and the partitioning objective of minimizing cutsize according to the connectivity metric definition given in (2.1). In this way, minimizing the cutsize corresponds to minimizing the upper bound given in Theorem 4.1 for the number of cache misses due to the access of $x$-*vector* entries. This reordering method will be referred to as "sHP$_{CN}$", where the small letter "s" is used to indicate the single-SpMxV framework.

**4.1.1. Exploiting temporal vs spatial locality in SpMxV.** Here we compare and contrast the existing HP-based method [43] (referred to as sHP$_{RN}$) and the proposed method sHP$_{CN}$ in terms exploiting temporal and spatial localities. Both sHP$_{RN}$ and sHP$_{CN}$ belong to the single-SpMxV framework and utilize 1D matrix partitioning for row/column reordering. For the CSR-based SpMxV operation, the row-net model utilized by sHP$_{RN}$ corresponds to the spatial locality hypergraph model proposed by Strout et al. [35] for data reordering of unstructured mesh computations. On the other hand, the column-net model utilized by sHP$_{CN}$ corresponds to the temporal locality hypergraph proposed by Strout et al. [35] for iteration reordering. Note that in the CSR-based SpMxV, the inner products of sparse rows with the dense input vector $x$ correspond to the iterations to be reordered. So the major difference between the sHP$_{RN}$ and sHP$_{CN}$ methods is that sHP$_{RN}$ primarily considers exploiting spatial locality and secondarily temporal locality, whereas sHP$_{CN}$ considers vice versa.

The above-mentioned difference between sHP$_{RN}$ and sHP$_{CN}$ can also be observed by investigating the row-net and column-net models used in these two HP-based methods. In HP with connectivity metric, the objective of cutsize minimization corresponds to clustering vertices with similar net connectivity to the same vertex parts. Hence, sHP$_{RN}$ clusters columns with similar sparsity patterns to the same column slice for partial column reordering thus exploiting spatial locality primarily, whereas sHP$_{CN}$ clusters rows with similar sparsity patterns to the same row slice for partial row reordering thus exploiting temporal locality primarily. In sHP$_{RN}$, the uncut and cut nets of a partition are used to decode the partial row reordering thus exploiting temporal locality secondarily. In sHP$_{CN}$, the uncut and cut nets of a partition are used to decode the partial column reordering thus exploiting spatial locality secondarily.

We should also note that the row-net and column-net models become equivalent for symmetric matrices. So, sHP$_{RN}$ and sHP$_{CN}$ obtain the same vertex partitions for symmetric matrices. The difference between these two methods in reordering matrices stems from the difference in the way that they decode the resultant partitions. sHP$_{RN}$ reorders the columns corresponding to the vertices in the same part of a partition successively, whereas sHP$_{CN}$ reorders the rows corresponding to the vertices in the same part of a partition successively.

**4.2. Row/column reordering based on 2D matrix partitioning.** We consider a row/column reordering which permutes a given matrix $A$ into a $K$-way DB form

$$\hat{A} = A_{DB} = P_r A P_c = \begin{bmatrix} A_{11} & & & & A_{1B} \\ & A_{22} & & & A_{2B} \\ & & \ddots & & \vdots \\ & & & A_{KK} & A_{KB} \\ A_{B1} & A_{B2} & \dots & A_{BK} & A_{BB} \end{bmatrix} = \begin{bmatrix} R_1 \\ R_2 \\ \vdots \\ R_K \\ R_B \end{bmatrix} = \begin{bmatrix} A'_{SB} \\ R_B \end{bmatrix}$$

$$= \begin{bmatrix} C_1 & C_2 & \dots & C_K & C_B \end{bmatrix}. \qquad (4.5)$$

Here, $R_B = [A_{B1}\, A_{B2}\, \dots\, A_{BK}\, A_{BB}]$ denotes the row border. Each row in $R_B$ is called a *column-coupling row* or simply a *coupling row*. $A'_{SB}$ denotes the columnwise SB part of $A_{DB}$ excluding the row border $R_B$. $R_k$ denotes the $k$th row slice of both $A'_{SB}$ and $A_{DB}$. $\lambda'(c_j)$ denotes the connectivity of column $c_j$ in $A'_{SB}$. $C'_B$ denotes the column border of $A'_{SB}$, whereas $C_B = [C'^T_B\, A^T_{BB}]^T$ denotes the column border of $A_{DB}$. $C_k = \begin{bmatrix} 0\dots 0\; A^T_{kk}\; 0\dots 0\; A^T_{Bk} \end{bmatrix}^T$ denotes the $k$th column slice of $A_{DB}$.

THEOREM 4.2. *Given a K-way DB form $A_{DB}$ of matrix $A$ such that every submatrix $R_k$ of $A'_{SB}$ fits into the cache, then we have*

$$\Phi_x(A_{DB}) \leq \sum_{c_j \in A'_{SB}} \lambda'(c_j) + \sum_{r_i \in R_B} nnz(r_i). \qquad (4.6)$$

*Proof.* We can consider the $y \leftarrow Ax$ multiply as two output-independent but input-dependent SpMxVs: $y_{SB} \leftarrow A'_{SB} x$ and $y_B \leftarrow R_B x$, where $y = [y^T_{SB}\, y^T_B]^T$. Thus $\Phi_x(A_{DB}) \leq \Phi_x(A'_{SB}) + \Phi_x(R_B)$. This upper bound occurs when no cache reuse occurs in accessing $x$-*vector* entries between the former and latter SpMxV operations. By proof of Theorem 4.1, we already have $\Phi_x(A'_{SB}) \leq \sum_{c_j} \lambda'(c_j)$. In the $y_B \leftarrow R_B x$ multiply, we have at most $nnz(r_i)$ $x$-vector accesses for each column-coupling row $r_i$ of $R_B$. This worst case happens when no cache-reuse occurs in accessing $x$-*vector* entries during the $y_B \leftarrow R_B x$ multiply. Hence, $\Phi_x(R_B) \leq \sum_{r_i \in R_B} nnz(r_i)$ thus concluding the proof. $\square$

Theorem 4.2 leads us to a cache-size-aware top-down row/column reordering through an $A$-to-$A_{DB}$ transformation that minimizes the upper bound given in (4.6) for $\Phi_x(A_{DB})$. Here, minimizing this sum relates to minimizing the number of cache misses due to the loss of temporal locality.

Here we propose to formulate the above-mentioned $A$-to-$A_{DB}$ transformation problem as an HP problem using the row-column-net model of matrix $A$ with a part size constraint of cache size. In the proposed formulation, column nets are associated with unit cost (i.e., $cost(n^c_j) = 1$ for each column $c_j$ and the cost of each row net is set to the number of nonzeros in the respective row (i.e., $cost(n^r_i) = nnz(r_i)$). However, existing HP tools do not handle a cutsize definition that encapsulates the right-hand side of (4.6), because the connectivity metric should be enforced for column nets, whereas the cut-net metric should be enforced for row nets. In order to encapsulate this different cutsize definition, we adapt and enhance the cut-net removal and cut-net splitting techniques adopted in RB algorithms utilized in HP tools. The connectivity of a column net should be calculated in such a way that it is as close as possible to the connectivity of the respective coupling column in the $A'_{SB}$ part of $A_{DB}$. For this purpose, after each bipartitioning step, each cut row-net is removed together with all vertices that it connects in both sides of the bipartition. Recall that the vertices of a cut net are not removed in the conventional cut-net removal scheme [9]. After applying the proposed removal scheme on the row nets on the cut, the conventional cut-net splitting technique [9] is applied to the column nets on the cut of the bipartition. This enhanced row-column-net model will be abbreviated as the "eRCN" model and the resulting reordering method will be referred to as "sHP$_{\text{eRCN}}$".

**Algorithm 3** SpMxV algorithm utilizing the multiple-SpMxV framework
---
**Require:** $A = A^1 + A^2 + \cdots + A^K$ partitioning of matrix $A$ and dense input vector $x$
**Output:** dense vector $y$
1: $y \leftarrow 0^T$
2: **for** $k \leftarrow 1$ **to** $K$ **do**
3: 　　$y \leftarrow y + A^k x$
4: **end for**
5: **return** y
---

The $K$-way partition $\{\mathcal{V}_1, \ldots, \mathcal{V}_K\}$ of $\mathcal{H}_{RCN}(A)$ obtained as a result of the abovementioned RB process is decoded as follows to induce a desired DB form of matrix $A$. The rows corresponding to the cut row-nets are permuted to the end to constitute the coupling rows of the row border $R_B$. The rows corresponding to the internal row-nets of part $\mathcal{V}_k$ are permuted to the $k$th row slice $R_k$. The columns corresponding to the internal column-nets of part $\mathcal{V}_k$ are permuted to the $k$th column slice $C_k$. It is clear that the columns corresponding to the cut column-nets remain in the column border $C_B$ of $A_{DB}$ and hence only those columns have the potential to remain in the column border $C'_B$ of $A'_{SB}$. Some of these columns may be permuted to a column slice $C_k$ if all of its nonzeros become confined to row slice $R_k$ and row border $R_B$. Such cases may occur as follows: Consider a cut column-net $n_j^c$ of a bipartition obtained at a particular RB step. If the internal row nets that belong to one part of the bipartition and that share a vertex with $n_j^c$ all become cut nets in the following RB steps, then column $c_j$ may be no longer a coupling column and be safely permuted to column slice $C_k$. For such cases, the proposed scheme fails to correctly encapsulate the column connectivity cost in $A'_{SB}$. The proposed cut row-net removal scheme avoids such column-connectivity miscalculations that may occur in the future RB steps due to the cut row-nets of the current bipartition. However, it is clear that our scheme cannot avoid such possible errors (related to the cut column-nets of the current bipartition) that may occur due to the row nets to be cut in the future RB steps.

**5. Multiple-SpMxV framework.** Let $\Pi = \{A^1, A^2, \ldots, A^K\}$ denote a splitting of matrix $A$ into $K$ $A^k$ matrices, where $A = A^1 + A^2 + \cdots + A^K$. In $\Pi$, $A^k$ matrices are mutually nonzero-disjoint, however they are not necessarily row disjoint or column disjoint. Note that every splitting $\Pi$ defines an access order on the matrix non-zeros, and every access order can define $\Pi$ that causes it.

In this framework, $y \leftarrow Ax$ operation is computed as a sequence of $K$ input- and output-dependent SpMxV operations, $y \leftarrow y + A^k x$ for $k = 1, \ldots, k$ as also shown in Algorithm 3. Individual SpMxV results are accumulated in the output vector $y$ on the fly in order to avoid additional write operations. The individual SpMxV operations are input dependent because of the shared columns among the $A^k$ matrices, whereas they are output dependent because of the shared rows among the $A^k$ matrices. Note that $A^k$ matrices are likely to contain empty rows and columns. The splitting of matrix $A$ should be done in such a way that the temporal and spatial localities of individual SpMxVs are exploited in order to minimize the number of cache misses.

We discuss pros and cons of this framework compared to the single-SpMxV framework in Section 5.1. In Section 5.2, we present a theorem that gives the guidelines for a "good" cache-size-aware matrix splitting based on 2D matrix partitioning. This theorem provides an upper bound on the total number of cache misses due to the access of $x$-$vector$ and $y$-$vector$ entries in all $y \leftarrow y + A^k x$ operations. The order of individual SpMxV operations is also important to exploit temporal locality between consecutive $y \leftarrow y + A^k x$ operations. In Section 5.3, we propose and discuss two methods for ordering SpMxV operations: RB ordering and TSP ordering.

**5.1. Pros and cons compared to single-SpMxV framework.** The single-SpMxV framework can be considered as a special case of multiple-SpMxV framework in which $A^k$ matrices are restricted to be row disjoint. Thus, the multiple-SpMxV framework brings an additional flexibility for exploiting the temporal and spatial localities. Clustering $A$-matrix rows/subrows with similar sparsity pattern into the same $A^k$ matrices increases the possibility of exploiting temporal locality in accessing $x$-$vector$ entries. Clustering $A$-matrix columns/subcolumns with similar sparsity pattern into the same $A^k$ matrices increases the possibility of exploiting spatial locality in accessing $x$-$vector$ entries as well as temporal locality in accessing $y$-$vector$ entries.

It is clear that single-SpMxV framework utilizing the CSR scheme severely suffers from dense rows. Dense rows causes loading large number of $x$-$vector$ entries to the cache thus disturbing the temporal locality in accessing $x$-$vector$ entries. The multiple-SpMxV framework may overcome this deficiency of the single-SpMxV framework through utilizing the flexibility of distributing the nonzeros of dense rows among multiple $A^k$ matrices in such a way to exploit the temporal locality in the respective $y \leftarrow y + A^k x$ operations.

However, this additional flexibility comes at a cost of disturbing the following localities compared to single SpMxV approach. There is some disturbance in the spatial locality in accessing the nonzeros of the $A$-matrix due to the division of three arrays associated with nonzeros into $K$ parts. However, this disturbance in spatial locality is negligible since elements of each of the three arrays are stored and accessed consecutively during each SpMxV operation. That is, at most $3(K-1)$ extra cache misses occur compared to the single SpMxV scheme due to the disturbance of spatial locality in accessing the nonzeros of $A$-matrix. More importantly, multiple read/writes of the individual SpMxV results might bring some disadvantages compared to single SpMxV scheme. These multiple read/writes disturb the spatial locality of accessing $y$-$vector$ entries as well as introducing a temporal locality exploitation problem in $y$-$vector$ entries.

**5.2. Splitting $A$ into $A^k$ matrices based on 2D matrix partitioning.** Given a splitting $\Pi$ of matrix $A$, let $\Phi_x(A, \Pi)$ and $\Phi_y(A, \Pi)$ respectively denote the number of cache misses due to the access of $x$-$vector$ and $y$-$vector$ entries during $y \leftarrow y + A^k x$ operations for $k = 1, \ldots, K$. Here, the total number of cache misses can be expressed as $\Phi(A, \Pi) = \Phi_x(A, \Pi) + \Phi_y(A, \Pi)$. Let $\lambda(r_i)$ and $\lambda(c_j)$ respectively denote the number of $A^k$ matrices that contain at least one nonzero of row $r_i$ and one nonzero of column $c_j$ of matrix $A$, i.e.,

$$\lambda(r_i) = |\{A^k \in \Pi : r_i \in A^k\}| \qquad (5.1a)$$
$$\lambda(c_j) = |\{A^k \in \Pi : c_j \in A^k\}|. \qquad (5.1b)$$

THEOREM 5.1. *Given a splitting* $\Pi = \{A^1, A^2, \ldots, A^K\}$ *of matrix* $A$, *then we have*

a) $\Phi_x(A, \Pi) \leq \sum_{c_j \in A} \lambda(c_j)$, *if each $A^k$ matrix fits into the cache;*

b) $\Phi_y(A, \Pi) \leq \sum_{r_i \in A} \lambda(r_i)$.

*Proof a).* Since each matrix $A^k$ fits into the cache, for any $c_j \in A^k$, the number of cache misses due to the access of $x_j$ is at most $\lambda(c_j)$ during all $y \leftarrow y + A^k x$ operations. This worst case happens when no cache-reuse occurs in accessing $x_j$ during successive $y \leftarrow y + A^k x$ operations.

*Proof b).* For any $r_i \in A^k$, the number of cache misses due to the access of $y_i$ is at most $\lambda(r_i)$ during all $y \leftarrow y + A^k x$ operations due to the nature of CSR-based SpMxV computation. This worst case happens when no cache-reuse occurs in accessing $y_i$ during successive $y \leftarrow y + A^k x$ operations. Note that the requirement of each $A^k$ fitting into the cache is not needed here. □

COROLLARY 5.2. *If each $A^k$ in $\Pi$ fits into the cache, then we have*

$$\Phi(A, \Pi) \leq \sum_{r_i \in A} \lambda(r_i) + \sum_{c_j \in A} \lambda(c_j). \tag{5.2}$$

Corollary 5.2 leads us to a cache-size-aware top-down matrix splitting that minimizes the upper bound given in (5.2) for $\Phi(A, \Pi)$. Here, minimizing this sum relates to minimizing the number of cache misses due to the loss of temporal locality.

The matrix splitting problem can be formulated as an HP-based 2D matrix parititioning using the row-column-net model [7, 11] of matrix $A$ with a part size constraint of cache size and partitioning objective of minimizing cutsize according to the connectivity metric definition given in (2.1). In this way, minimizing the cutsize corresponds to minimizing the upper bound given in Theorem 5.1 for the total number of cache misses due to the access of $x$-*vector* and $y$-*vector* entries. This reordering method will be referred to as "mHP$_{\text{RCN}}$", where the small letter "m" is used to indicate the multiple-SpMxV framework.

**5.3. Ordering individual SpMxV operations.** The above-mentioned objective in splitting matrix $A$ into $A^k$ matrices is to exploit temporal locality of individual SpMxVs in order to minimize the number of cache misses. However, when all SpMxVs are considered, data reuse between two consecutive SpMxVs should better be considered to exploit temporal locality. Here we propose and discuss two methods for ordering SpMxV operations: RB ordering and TSP ordering.

**5.3.1. RB ordering.** The RB tree constructed during the recursive hypergraph bipartitioning is a full binary tree, where each node represents a vertex subset as well as the respective induced sub-hypergraph on which a 2-way HP to be applied. Note that the root node represents both the original vertex set and the original row-column-net hypergraph model for the given $A$-matrix and the leaf nodes represent the $A^k$ matrices. The preorder, postorder and inorder traversals starting from the root node give the same traversal order on the leaf nodes, thus inducing an RB order on the individual SpMxV operations of the multiple-SpMxV framework. In the RB tree, the amount of row/column sharing between two leaf nodes ($A^k$ matrices) is expected to decrease with increasing path length to their first common ancestor in the RB tree. Since sibling nodes have a common parent, the $A^k$ matrices corresponding to the sibling leaf-node pairs are likely to share larger number of rows and columns compared to $A^k$ matrices corresponding to the non-sibling leaf node pairs. As this scheme orders the sibling leaf nodes consecutively, the RB ordering is expected to yield an order on $A^k$ matrices that respects temporal locality in accessing $x$-*vector* and $y$-*vector* entries.

**5.3.2. TSP ordering.** Let $\widehat{\Pi} = <A^1, A^2, \ldots, A^K>$ denote an ordered version of a given splitting $\Pi$. A subchain of $\widehat{\Pi}$ is said to cover a row $r_i$ and a column $c_j$ if each $A^k$ matrix in the subchain contains at least one nonzero of row $r_i$ and column $c_j$, respectively. Let $\gamma(r_i)$ and $\gamma(c_j)$ denote the number of maximal $A^k$-matrix subchains that cover row $r_i$ and column $c_j$, respectively. Let $L$ denote the cache line size. Let $\Phi(A, \widehat{\Pi})$ denote the total number of cache misses due to the access of $x$-*vector* and $y$-*vector* entries for a given order $\widehat{\Pi}$ of $y \leftarrow y + A^k x$ operations for $k = 1, \ldots, K$. Theorem 5.3 gives a lower bound for $\Phi(A, \widehat{\Pi})$ and Theorem 5.4 shows our TSP formulation that minimizes this lower bound.

THEOREM 5.3. *Given an ordered splitting $\widehat{\Pi}$ of matrix $A$, if none of the $A^k$ matrices fit into the cache, then we have*

$$\Phi(A, \widehat{\Pi}) \geq \frac{\sum_{r_i \in A} \gamma(r_i) + \sum_{c_j \in A} \gamma(c_j)}{L}. \quad (5.3)$$

*Proof.* We first consider the case $L = 1$. Consider a column $c_j$ of matrix $A$. Then there exists $\gamma(c_j)$ maximal $A^k$-matrix subchains that cover column $c_j$. Since none of the $A^k$ matrices can fit into the cache, it is guaranteed that there will be no cache reuse of column $c_j$ between two different maximal $A^k$-matrix subchains that cover $c_j$. Therefore, at least $\gamma(c_j)$ cache misses will occur for each column $c_j$ which means that $\Phi_x(A, \widehat{\Pi}) \geq \sum_{c_j} \gamma(c_j)$. A similar proof follows for a row $r_i$ of matrix $A$ so that $\Phi_y(A, \widehat{\Pi}) \geq \sum_{r_i} \gamma(r_i)$. When $L > 1$, the number of cache misses may decrease $L$-fold at most. □

We define the TSP instance $(\mathcal{G} = (\mathcal{V}, \mathcal{E}), w)$ over a given unordered splitting $\Pi$ of matrix $A$ as follows. The vertex set $\mathcal{V}$ denotes the set of $A^k$ matrices. The weight $w(k, \ell)$ of edge $e_{k\ell} \in \mathcal{E}$ is set to be equal to the sum of the number of shared rows and columns between $A^k$ and $A^\ell$.

THEOREM 5.4. *For a given unordered splitting $\Pi$ of matrix $A$, finding an order on the vertices of the TSP instance $(\mathcal{G}, w)$ that maximizes the path weight corresponds to finding an order $\widehat{\Pi}$ of $A^k$ matrices that minimizes the lower bound given in (5.4) for $\Phi(A, \widehat{\Pi})$.*

*Proof.* Below, let $A^{\Gamma(\ell)}$ denote the $\ell$th $A^k$ matrix in the ordering $\Gamma$ of $A^k$ matrices.

$$\Psi = \sum_{r_i}\left[|A^{\Gamma(1)} \cap \{r_i\}| + \sum_{\ell=2}^{K}|(A^{\Gamma(\ell)} - A^{\Gamma(\ell-1)}) \cap \{r_i\}|\right]$$
$$+ \sum_{c_j}\left[|A^{\Gamma(1)} \cap \{c_j\}| + \sum_{\ell=2}^{K}|(A^{\Gamma(\ell)} - A^{\Gamma(\ell-1)}) \cap \{c_j\}|\right]$$
$$= |A^{\Gamma(1)}| + \sum_{\ell=2}^{K}|(A^{\Gamma(\ell)} - A^{\Gamma(\ell-1)})| \quad = \quad |A^{\Gamma(1)}| + \sum_{\ell=2}^{K}(|A^{\Gamma(\ell)}| - |A^{\Gamma(\ell)} \cap A^{\Gamma(\ell-1)}|)$$
$$= \sum_{\ell=1}^{K}|A^{\Gamma(\ell)}| - \sum_{\ell=2}^{K}|A^{\Gamma(\ell)} \cap A^{\Gamma(\ell-1)}| \quad = \quad \sum_{\ell=1}^{K}|A^{\Gamma(\ell)}| - \sum_{\ell=2}^{K}w(\Gamma(\ell), \Gamma(\ell-1))$$

The maximum value of $\sum_{\ell=2}^{K} w(\Gamma(\ell), \Gamma(\ell-1))$ will yield the minimum value of $\sum_{r_i \in A} \gamma(r_i) + \sum_{c_j \in A} \gamma(c_j)$. Then, finding an order on $\mathcal{V}$ that maximizes the path weight $\sum_{\ell=2}^{K} w(\Gamma(\ell), \Gamma(\ell-1))$ corresponds to finding an order of submatrices that minimizes $\sum_{r_i \in A} \gamma(r_i) + \sum_{c_j \in A} \gamma(c_j)$. □

## 6. Experimental results.

**6.1. Experimental setup.** We tested the performance of the proposed methods against three state-of-the-art methods: sBFS [35], sRCM [13, 24, 38] and sHP$_{\text{RN}}$ [43] all of which belong to the single-SpMxV framework. Here, sBFS refers to our adaptation of BFS-based simultaneous data and iteration reordering method of Strout et al. [35] to matrix row and column reordering. Strout et al.'s method depends on implementing breadth-first search on both temporal and spatial locality hypergraphs simultaneously. In our adaptation, we apply BFS on the bipartite graph representation of the matrix, so that the resulting BFS orders on the row and column vertices determine row and column reorderings, respectively. sRCM refers to applying the RCM method, which is widely used for envelope reduction of symmetric matrices, on the bipartite graph representation of the given sparse matrix. Application of the RCM method to bipartite graphs has also been used by Berry et al. [4] to reorder rectangular term-by-document matrices for envelope minimization. sHP$_{RN}$ refers to the work by Yzelman and Bisseling [43] which utilizes HP using the row-net model for CSR-based SpMxV.

The HP-based top-down reordering methods sHP$_{RN}$, sHP$_{CN}$, sHP$_{eRCN}$ and mHP$_{RCN}$ are implemented using the state-of-the-art HP tool PaToH [10]. In these implementations, PaToH is used as a 2-way HP tool within the RB paradigm. The hypergraphs representing sparse matrices according to the respective models are recursively bipartitioned into parts until the CSR-storage size of the submatrix (together with the $x$ and $y$ vectors) corresponding to a part drops below the cache size. PaToH is used with default parameters except the use of heavy connectivity clustering (`PATOH_CRS_HCC=9`) in the sHP$_{RN}$, sHP$_{CN}$ and sHP$_{eRCN}$ methods that belong to the single-SpMxV framework, and the use of absorption clustering using nets (`PATOH_CRS_ABSHCC=11`) in the mHP$_{RCN}$ method that belong to the multiple-SpMxV framework. Since PaToH contains randomized algorithms, the reordering results are reported by averaging the values obtained in 10 different runs, each randomly seeded.

Performance evaluations are carried out in two different settings: cache-miss simulations and actual running times by using OSKI (BeBOP Optimized Sparse Kernel Interface Library) [41]. In cache-miss simulations, 8-byte words are used for matrix nonzeros, $x$-$vector$ entries and $y$-$vector$ entries. In OSKI runs, double precision arithmetic is used. Cache-miss simulations are performed on 36 small-to-medium size matrices, whereas OSKI runs are performed on 17 large size matrices. All test matrices are obtained from the University of Florida Sparse Matrix Collection [37]. CSR-storage sizes of small-to-medium size matrices vary between 441 KB to 13 MB, whereas CSR-storage sizes of large size matrices vary between 13 MB to 94 MB. Properties of these matrices are presented in Table 6.1. As seen in the table, both sets of small-to-medium and large size matrices are categorized into three groups as symmetric, square nonsymmetric and rectangular. In each group, the test matrices are listed in the order of increasing number of nonzeros ("nnz"). In the table, "avg" and "max" denote the average and the maximum number of nonzeros per row/column. "cov" denotes the coefficient of variation of number of nonzeros per row/column. The "cov" value of a matrix can be considered as an indication of the level of irregularity in the number of nonzeros per row and column.

**6.2. Cache-miss simulations.** Cache-miss simulations are performed by running the single-level cache simulator developed by Yzelman and Bisseling [43] on small-to-medium size test matrices. The simulator is configured to have 64 KB, 2-way set-associative cache with a line size of 64 bytes (8 words). Some of the experiments are conducted to show the sensitivities of the methods to the cache-line size without changing the other cache parameters. In the simulations, since the ICSR [27] storage scheme is to be used in the multiple-SpMxV framework as discussed in Section 5, ICSR is also used for all other methods. The ZZCSR scheme proposed by Yzelman and Bisseling [43] is not used in the simulations, since the main purpose of this work is to show the cache miss effects of the six different reordering methods. In Tables 6.3, 6.4, 6.5 and 6.8, the performances of the existing and proposed methods are displayed in terms of normalized cache miss values, where each normalized value is calculated through dividing the number of cache misses for the reordered matrix by that of the original matrix. In these tables, the "$x$", "$y$" and "$x+y$" columns, respectively, denote the normalized cache miss values due to the access of $x$-$vector$ entries, $y$-$vector$ entries and both. In these tables, compulsory cache misses due to the access of matrix nonzeros are not reported in order to better show the performance differences among the methods.

TABLE 6.1
*Properties of test matrices*

| Name | number of | | | nnz's in a row | | | nnz's in a column | | |
|---|---|---|---|---|---|---|---|---|---|
| | rows | cols | nonzeros | avg | max | cov | avg | max | cov |
| *Small-to-Medium Size Matrices* | | | | | | | | | |
| *Symmetric Matrices* | | | | | | | | | |
| ncvxqp9 | 16,554 | 16,554 | 54,040 | 3 | 9 | 0.5 | 3 | 9 | 0.5 |
| tuma1 | 22,967 | 22,967 | 87,760 | 4 | 5 | 0.3 | 4 | 5 | 0.3 |
| bloweybl | 30,003 | 30,003 | 120,000 | 4 | 10,001 | 14.4 | 4 | 10,001 | 14.4 |
| bloweya | 30,004 | 30,004 | 150,009 | 5 | 10,001 | 11.6 | 5 | 10,001 | 11.6 |
| brainpc2 | 27,607 | 27,607 | 179,395 | 7 | 13,799 | 20.2 | 7 | 13,799 | 20.2 |
| a5esindl | 60,008 | 60,008 | 255,004 | 4 | 9,993 | 12.7 | 4 | 9,993 | 12.7 |
| dixmaanl | 60,000 | 60,000 | 299,998 | 5 | 5 | 0.0 | 5 | 5 | 0.0 |
| shallow_water1 | 81,920 | 81,920 | 327,680 | 4 | 4 | 0.0 | 4 | 4 | 0.0 |
| c-65 | 48,066 | 48,066 | 360,528 | 8 | 3,276 | 2.5 | 8 | 3,276 | 2.5 |
| finan512 | 74,752 | 74,752 | 596,992 | 8 | 55 | 0.8 | 8 | 55 | 0.8 |
| copter2 | 55,476 | 55,476 | 759,952 | 14 | 45 | 0.3 | 14 | 45 | 0.3 |
| msc23052 | 23,052 | 23,052 | 1,154,814 | 50 | 178 | 0.2 | 50 | 178 | 0.2 |
| *Square Nonsymmetric Matrices* | | | | | | | | | |
| poli_large | 15,575 | 15,575 | 33,074 | 2 | 491 | 4.2 | 2 | 18 | 0.2 |
| powersim | 15,838 | 15,838 | 67,562 | 4 | 40 | 0.6 | 4 | 41 | 0.8 |
| memplus | 17,758 | 17,758 | 126,150 | 7 | 574 | 3.1 | 7 | 574 | 3.1 |
| Zhao1 | 33,861 | 33,861 | 166,453 | 5 | 6 | 0.1 | 5 | 7 | 0.2 |
| mult_dcop_01 | 25,187 | 25,187 | 193,276 | 8 | 22,767 | 18.7 | 8 | 22,774 | 18.8 |
| jan99jac120sc | 41,374 | 41,374 | 260,202 | 6 | 68 | 1.1 | 6 | 138 | 2.3 |
| circuit_4 | 80,209 | 80,209 | 307,604 | 4 | 6,750 | 7.8 | 4 | 8,900 | 10.5 |
| ckt11752_dc_1 | 49,702 | 49,702 | 333,029 | 7 | 2,921 | 3.5 | 7 | 2,921 | 3.5 |
| poisson3Da | 13,514 | 13,514 | 352,762 | 26 | 110 | 0.5 | 26 | 110 | 0.5 |
| bcircuit | 68,902 | 68,902 | 375,558 | 6 | 34 | 0.4 | 6 | 34 | 0.4 |
| g7jac120 | 35,550 | 35,550 | 475,296 | 13 | 153 | 1.7 | 13 | 120 | 1.7 |
| e40r0100 | 17,281 | 17,281 | 553,562 | 32 | 62 | 0.5 | 32 | 62 | 0.5 |
| *Rectangular Matrices* | | | | | | | | | |
| lp_dfl001 | 6,071 | 12,230 | 35,632 | 6 | 228 | 1.3 | 3 | 14 | 0.4 |
| ge | 10,099 | 16,369 | 44,825 | 4 | 48 | 0.8 | 3 | 36 | 0.9 |
| ex3sta1 | 17,443 | 17,516 | 68,779 | 4 | 8 | 0.4 | 4 | 46 | 1.4 |
| lp_stocfor3 | 16,675 | 23,541 | 76,473 | 5 | 15 | 0.7 | 3 | 18 | 1.0 |
| cq9 | 9,278 | 21,534 | 96,653 | 10 | 391 | 3.5 | 5 | 24 | 1.0 |
| psse0 | 26,722 | 11,028 | 102,432 | 4 | 4 | 0.1 | 9 | 68 | 0.7 |
| co9 | 10,789 | 22,924 | 109,651 | 10 | 441 | 3.6 | 5 | 28 | 1.1 |
| baxter | 27,441 | 30,733 | 111,576 | 4 | 2,951 | 8.7 | 4 | 46 | 1.6 |
| graphics | 29,493 | 11,822 | 117,954 | 4 | 4 | 0.0 | 10 | 87 | 1.0 |
| fome12 | 24,284 | 48,920 | 142,528 | 6 | 228 | 1.3 | 3 | 14 | 0.4 |
| route | 20,894 | 43,019 | 206,782 | 10 | 2,781 | 7.1 | 5 | 44 | 1.0 |
| fxm4_6 | 22,400 | 47,185 | 265,442 | 12 | 57 | 1.0 | 6 | 24 | 1.1 |
| *Large Size Matrices* | | | | | | | | | |
| *Symmetric Matrices* | | | | | | | | | |
| c-73 | 169,422 | 169,422 | 1,279,274 | 8 | 39,937 | 20.1 | 8 | 39,937 | 20.1 |
| c-73b | 169,422 | 169,422 | 1,279,274 | 8 | 39,937 | 20.1 | 8 | 39,937 | 20.1 |
| rgg_n_2_17_s0 | 131,072 | 131,072 | 1,457,506 | 11 | 96 | 0.3 | 11 | 28 | 0.3 |
| boyd2 | 466,316 | 466,316 | 1,500,397 | 3 | 93,262 | 60.6 | 3 | 93,262 | 60.6 |
| ins2 | 309,412 | 309,412 | 2,751,484 | 9 | 303,879 | 65.3 | 9 | 309,412 | 66.4 |
| rgg_n_2_18_s0 | 262,144 | 262,144 | 3,094,566 | 12 | 62 | 0.3 | 12 | 31 | 0.3 |
| *Square Nonsymmetric Matrices* | | | | | | | | | |
| Raj1 | 263,743 | 263,743 | 1,302,464 | 5 | 40,468 | 17.9 | 5 | 40,468 | 17.9 |
| rajat21 | 411,676 | 411,676 | 1,893,370 | 5 | 118,689 | 41.0 | 5 | 100,470 | 34.8 |
| rajat24 | 358,172 | 358,172 | 1,948,235 | 5 | 105,296 | 33.1 | 5 | 105,296 | 33.1 |
| ASIC_320k | 321,821 | 321,821 | 2,635,364 | 8 | 203,800 | 61.4 | 8 | 203,800 | 61.4 |
| Stanford_Berkeley | 683,446 | 683,446 | 7,583,376 | 11 | 76,162 | 25.0 | 11 | 249 | 1.5 |
| *Rectangular Matrices* | | | | | | | | | |
| kneser_10_4_1 | 349,651 | 330,751 | 992,252 | 3 | 51,751 | 31.9 | 3 | 3 | 0.0 |
| neos | 479,119 | 515,905 | 1,526,794 | 3 | 29 | 0.2 | 3 | 16,220 | 15.6 |
| wheel_601 | 902,103 | 723,605 | 2,170,814 | 2 | 442,477 | 193.9 | 3 | 3 | 0.0 |
| LargeRegFile | 2,111,154 | 801,374 | 4,944,201 | 2 | 4 | 0.3 | 6 | 655,876 | 145.9 |
| cont1_l | 1,918,399 | 1,921,596 | 7,031,999 | 4 | 5 | 0.3 | 4 | 1,279,998 | 252.3 |
| degme | 185,501 | 659,415 | 8,127,528 | 44 | 624,079 | 33.1 | 12 | 18 | 0.1 |

Table 6.2 displays the empirical corroboration for the proposed theorems on large size matrices. The values given in the table are obtained by running cache simulations for a 2 MB, 8-way set-associative cache with a line size of 8 bytes. In this table, second, third and fourth columns show the numbers of cache misses $\widetilde{\Phi}_x(A_{SB})$, $\widetilde{\Phi}_x(A_{DB})$ and $\widetilde{\Phi}(A, \Pi)$ of the sHP$_{\text{CN}}$, sHP$_{\text{eRCN}}$ and mHP$_{\text{RCN}}$ methods normalized with respect to the upper bounds given in Theorems 4.1, 4.2 and 5.1, respectively. The last column of this table shows the number of cache misses $\widetilde{\Phi}(A, \widehat{\Pi})$ of the mHP$_{\text{RCN}}$ method normalized with respect to the lower bound given in Theorem 5.3. The values close to 1.00 for $\widetilde{\Phi}_x(A_{SB})$, $\widetilde{\Phi}(A, \Pi)$ and $\widetilde{\Phi}(A, \widehat{\Pi})$ show that the bounds given in Theorems 4.1, 5.1 and 5.3 are quite tight. For mHP$_{\text{RCN}}$, a few $\widetilde{\Phi}(A, \Pi)$ values slightly larger than 1.00 can be attributed to the fact that the cache is not fully associative in the simulations. Rather small values for $\widetilde{\Phi}_x(A_{DB})$, show that the upper bound given in Theorem 4.2 is rather loose. This can be attributed to uncontrolled temporal locality in $y_B \leftarrow R_B x$ due to the possibility of rows with similar sparsity pattern in $R_B$.

TABLE 6.2
*Empirical corroboration for Theorems 4.1, 4.2, 5.1 and 5.3*

|  | Proposed Methods | | | |
|---|---|---|---|---|
|  | Single SpMxV | | Mult. SpMxVs | |
|  | sHP$_{\text{CN}}$ | sHP$_{\text{eRCN}}$ | mHP$_{\text{RCN}}$ | |
|  | $\widetilde{\Phi}_x(A_{SB})$ | $\widetilde{\Phi}_x(A_{DB})$ | $\widetilde{\Phi}(A, \Pi)$ | $\widetilde{\Phi}(A, \widehat{\Pi})$ |
| Symmetric Matrices | | | | |
| c-73 | 0.85 | 0.35 | 0.99 | 1.00 |
| c-73b | 0.85 | 0.35 | 0.99 | 1.00 |
| rgg_n_2_17_s0 | 0.99 | 0.49 | 0.99 | 1.00 |
| boyd2 | 0.98 | 0.58 | 1.00 | 1.00 |
| ins2 | 0.96 | 0.71 | 1.00 | 1.00 |
| rgg_n_2_18_s0 | 0.99 | 0.49 | 0.99 | 1.00 |
| Square Nonsymmetric Matrices | | | | |
| Raj1 | 0.92 | 0.55 | 1.00 | 1.00 |
| rajat21 | 0.97 | 0.58 | 1.02 | 1.02 |
| rajat24 | 0.95 | 0.58 | 1.01 | 1.01 |
| ASIC_320k | 0.96 | 0.67 | 1.00 | 1.00 |
| Stanford_Berkeley | 0.69 | 0.30 | 0.95 | 0.95 |
| Rectangular Matrices | | | | |
| kneser_10_4_1 | 0.99 | 0.53 | 0.99 | 1.00 |
| neos | 0.99 | 0.50 | 1.00 | 1.00 |
| wheel_601 | 0.99 | 0.59 | 1.00 | 1.00 |
| LargeRegFile | 1.00 | 0.50 | 1.00 | 1.00 |
| cont1_l | 1.00 | 0.50 | 1.00 | 1.00 |
| degme | 0.99 | 0.65 | 1.02 | 1.02 |
| Geometric Means | | | | |
| Symmetric | 0.97 | 0.49 | 0.99 | 1.00 |
| Nonsymmetric | 0.94 | 0.56 | 1.00 | 1.00 |
| Rectangular | 0.99 | 0.52 | 1.00 | 1.00 |
| Overall | 0.97 | 0.52 | 1.00 | 1.00 |

We introduce Table 6.3 to show the validity of the enhanced row-column-net model proposed in Section 4.2 for the sHP$_{\text{eRCN}}$ method. In the table, sHP$_{\text{RCN}}$ refers to a version of the sHP$_{\text{eRCN}}$ method that utilizes the conventional row-column-net model instead of the enhanced row-column-net model. Table 6.3 displays average performance results of sHP$_{\text{RCN}}$ and sHP$_{\text{eRCN}}$ over the three different matrix categories as well as the overall averages. As seen in the table, sHP$_{\text{eRCN}}$ performs considerably better than sHP$_{\text{RCN}}$, thus showing the validity of the cutsize definition that encapsulates the right-hand side of (4.6).

CACHE LOCALITY IN SPARSE-MATRIX VECTOR MULTIPLY 19TABLE 6.3
*Average simulation results (misses) to display the merits of enhancement of the row-column-net model in* $sHP_{eRCN}$ *(cache size = part-weight threshold = 64 KB)*

|  | $sHP_{RCN}$ | $sHP_{eRCN}$ |
|---|---|---|
|  | $x$ | $x$ |
| Symmetric | 0.54 | 0.47 |
| Nonsymmetric | 0.45 | 0.40 |
| Rectangular | 0.44 | 0.43 |
| Overall | 0.48 | 0.43 |

TABLE 6.4
*Average simulation results (misses) to display the merits of ordering SpMxV operations in* $mHP_{RCN}$ *(cache size = part-weight threshold = 64 KB)*

|  | Random Ordering | | | RB Ordering | | | TSP Ordering | | |
|---|---|---|---|---|---|---|---|---|---|
|  | $x$ | $y$ | $x+y$ | $x$ | $y$ | $x+y$ | $x$ | $y$ | $x+y$ |
| Symmetric | 0.44 | 1.34 | 0.62 | 0.41 | 1.28 | 0.58 | 0.40 | 1.26 | 0.57 |
| Nonsymmetric | 0.37 | 1.60 | 0.54 | 0.34 | 1.55 | 0.50 | 0.34 | 1.54 | 0.50 |
| Rectangular | 0.27 | 1.39 | 0.40 | 0.26 | 1.35 | 0.39 | 0.27 | 1.36 | 0.40 |
| Overall | 0.35 | 1.44 | 0.51 | 0.33 | 1.39 | 0.49 | 0.33 | 1.38 | 0.48 |

We introduce Table 6.4 to show the merits of ordering individual SpMxV operations in the $mHP_{RCN}$ method. The table displays average performance results of $mHP_{RCN}$ for the random, RB and TSP orderings over the three different matrix categories as well as the overall averages. As seen in the table, both RB and TSP orderings lead to considerable performance improvement in the $mHP_{RCN}$ method compared to the random ordering, where the TSP ordering leads to slightly better improvement than the RB ordering. In the following tables, we display the performance results of the $mHP_{RCN}$ method that utilizes the TSP ordering. The TSP implementation given in [21] is used in these experiments.

Table 6.5 displays the performance comparison of the existing and proposed methods for small-to-medium size matrices. The bottom part of the table shows the geometric means of the performance results of the methods over the three different matrix categories as well as the overall averages. Among the existing methods, $sHP_{RN}$ performs considerably better than both sBFS and sRCM for all matrix categories, on the average.

**6.2.1. Comparison of 1D methods $sHP_{RN}$ and $sHP_{CN}$.** Here we present the experimental comparison of $sHP_{RN}$ and $sHP_{CN}$ and show how this experimental comparison relates to the theoretical comparison given in Section 4.1.1. As seen Table 6.5, $sHP_{CN}$ performs significantly better than $sHP_{RN}$, on the overall average. $sHP_{CN}$ performs better than $sHP_{RN}$ in all of the 36 reordering instances except `a5esindl`, `lp_stocfactor3` and `route`. The significant performance gap between $sHP_{RN}$ and $sHP_{CN}$ in favor of $sHP_{CN}$ even for symmetric matrices confirm our expectation that temporal locality is more important than spatial locality in SpMxV operations that involve irregularly sparse matrices.

We introduce Table 6.6 to experimentally investigate the sensitivity of the $sHP_{RN}$ and $sHP_{CN}$ methods to the cache-line size. In the construction of the averages reported in this table, simulation results of every method are normalized with respect to those of the original ordering with the respective cache-line size. We also utilize Table 6.6 to provide fairness in the comparison of $sHP_{RN}$ and $sHP_{CN}$ methods for nonsymmetric square and rectangular matrices. Some of the nonsymmetric square and rectangular matrices may be more suitable for rowwise partitioning by the column-net model, whereas some other matrices may be more suitable for columnwise partitioning utilizing the row-net model. This is because of the differences in row and column sparsity patterns of a given nonsymmetric or rectangular matrix. Hendrickson and Kolda [22] and Ucar and Aykanat [39] provide discussions on choosing partitioning dimension depending on the individual matrix characteristics in the parallel SpMxV context. In the construction of Table 6.6, each of the $sHP_{RN}$ and $sHP_{CN}$ methods are applied on both $A$ and $A^T$ matrices and the better result is reported for the respective method on the reordering of matrix $A$. Here the performance of CSR-based SpMxV $y \leftarrow A^T x$ is assumed to simulate the performance of CSC-based $y \leftarrow Ax$. Comparison of the results in Table 6.6 for the line size of 64 bytes and the average results in Table 6.5 shows that the performance of both methods increase due to the selection of better partitioning dimension (especially for rectangular matrices) while the performance gap remaining almost the same.

TABLE 6.5
*Simulation results (misses) for small-to-medium size test matrices (cache size = part-weight threshold = 64 KB)*

| | Existing Methods | | | | | | Proposed Methods | | | | | |
|---|---|---|---|---|---|---|---|---|---|---|---|---|
| | Single SpMxV | | | | | | | | | | Multiple SpMxVs | |
| | sBFS [35] | | sRCM [24] Modified | | sHP$_{RN}$ [43] (1D Part.) | | sHP$_{CN}$ (1D Part.) | | sHP$_{eRCN}$ (2D Part.) | | mHP$_{RCN}$ (2D Partioning) | | |
| | $x$ | $x+y$ | $x$ | $x+y$ | $x$ | $x+y$ | $x$ | $x+y$ | $x$ | $x+y$ | $x$ | $y$ | $x+y$ |
| | | | | | Symmetric Matrices | | | | | | | | |
| ncvxqp9 | 0.51 | 0.59 | 0.52 | 0.60 | 0.37 | 0.48 | 0.28 | **0.40** | 0.28 | **0.40** | 0.31 | 1.19 | 0.47 |
| tuma1 | 0.42 | **0.59** | 0.58 | 0.71 | 0.62 | 0.73 | 0.56 | 0.69 | 0.56 | 0.69 | 0.45 | 1.01 | 0.60 |
| bloweybl | 1.00 | 1.00 | 1.00 | 1.00 | 0.88 | 0.92 | 0.68 | 0.77 | 0.63 | **0.74** | 0.63 | 1.02 | **0.74** |
| bloweya | 1.00 | 1.00 | 1.03 | 1.02 | 1.18 | 1.12 | 0.65 | 0.75 | 0.73 | 0.81 | 0.45 | 1.03 | **0.62** |
| brainpc2 | 0.88 | 0.90 | 0.89 | 0.91 | 1.33 | 1.27 | 1.08 | 1.06 | 0.66 | 0.73 | 0.28 | 1.04 | **0.43** |
| a5esindl | 1.11 | 1.09 | 0.83 | 0.86 | 0.84 | 0.87 | 1.12 | 1.10 | 0.40 | 0.52 | 0.28 | 1.03 | **0.43** |
| dixmaanl | 0.33 | **0.50** | 0.33 | **0.50** | 0.34 | 0.51 | 0.34 | **0.50** | 0.34 | **0.50** | 0.36 | 1.01 | 0.52 |
| shallow_water1 | 1.45 | 1.28 | 1.22 | 1.14 | 1.10 | 1.07 | 0.90 | 0.94 | 0.89 | 0.94 | 0.77 | 1.01 | **0.86** |
| c-65 | 0.90 | 0.91 | 0.96 | 0.97 | 0.61 | 0.67 | 0.38 | 0.47 | 0.35 | 0.44 | 0.24 | 1.37 | **0.40** |
| finan512 | 1.57 | 1.40 | 1.47 | 1.34 | 0.65 | 0.75 | 0.56 | **0.68** | 0.55 | **0.68** | 0.70 | 1.27 | 0.89 |
| copter2 | 0.44 | 0.49 | 0.43 | 0.49 | 0.41 | 0.47 | 0.26 | **0.33** | 0.26 | **0.33** | 0.30 | 2.60 | 0.53 |
| msc23052 | 0.46 | 0.51 | 0.42 | 0.47 | 0.52 | 0.57 | 0.40 | **0.46** | 0.44 | 0.49 | 0.35 | 2.65 | 0.57 |
| | | | | | Square Nonsymmetric Matrices | | | | | | | | |
| poli_large | 1.12 | 1.08 | 1.12 | 1.08 | 0.86 | 0.91 | 0.62 | **0.75** | 0.64 | 0.77 | 0.60 | 1.05 | 0.76 |
| powersim | 1.02 | 1.01 | 0.72 | 0.81 | 0.55 | 0.69 | 0.51 | **0.66** | 0.51 | **0.66** | 0.50 | 1.04 | 0.67 |
| memplus | 0.87 | 0.90 | 1.06 | 1.05 | 1.39 | 1.30 | 0.91 | 0.93 | 0.87 | 0.90 | 0.47 | 1.24 | **0.63** |
| Zhao1 | 0.55 | 0.65 | 0.36 | **0.51** | 0.72 | 0.79 | 0.48 | 0.60 | 0.49 | 0.60 | 0.60 | 1.64 | 0.84 |
| mult_dcop_01 | 0.98 | 0.98 | 1.00 | 1.00 | 0.70 | 0.71 | 0.45 | 0.48 | 0.18 | 0.23 | 0.13 | 1.36 | **0.20** |
| jan99jac120sc | 1.20 | 1.15 | 1.30 | 1.22 | 0.92 | 0.94 | 0.51 | **0.62** | 0.52 | 0.63 | 0.63 | 1.41 | 0.83 |
| circuit_4 | 1.52 | 1.39 | 1.83 | 1.62 | 1.45 | 1.34 | 0.94 | 0.95 | 0.87 | 0.91 | 0.41 | 1.16 | **0.60** |
| ckt11752_dc_1 | 0.79 | 0.83 | 0.86 | 0.89 | 0.58 | 0.66 | 0.40 | 0.52 | 0.42 | 0.54 | 0.31 | 1.10 | **0.47** |
| poisson3Da | 0.09 | 0.11 | 0.09 | 0.11 | 0.14 | 0.15 | 0.09 | **0.10** | 0.09 | **0.10** | 0.09 | 6.32 | 0.18 |
| bcircuit | 0.60 | 0.67 | 0.58 | 0.66 | 0.32 | 0.44 | 0.26 | **0.39** | 0.26 | **0.39** | 0.27 | 1.10 | 0.42 |
| g7jac120 | 0.75 | 0.76 | 0.26 | 0.30 | 0.44 | 0.47 | 0.21 | **0.25** | 0.23 | 0.28 | 0.18 | 2.51 | 0.31 |
| e40r0100 | 0.82 | 0.86 | 0.74 | 0.79 | 0.76 | 0.81 | 0.63 | **0.71** | 0.66 | 0.73 | 0.54 | 1.94 | 0.84 |
| | | | | | Rectangular Matrices | | | | | | | | |
| lp_dfl001 | 0.30 | 0.33 | 0.31 | 0.34 | 0.34 | 0.36 | 0.18 | 0.21 | 0.20 | 0.23 | 0.10 | 2.57 | **0.19** |
| ge | 0.40 | 0.47 | 0.34 | 0.41 | 0.30 | 0.37 | 0.25 | 0.33 | 0.25 | 0.33 | 0.21 | 1.23 | **0.32** |
| ex3sta1 | 1.75 | 1.47 | 1.14 | 1.09 | 1.23 | 1.14 | 0.86 | 0.91 | 0.81 | **0.88** | 0.81 | 1.09 | 0.91 |
| lp_stocfor3 | 1.74 | 1.48 | 1.64 | 1.42 | 0.79 | **0.86** | 0.80 | 0.87 | 0.80 | 0.87 | 0.79 | 1.02 | 0.87 |
| cq9 | 0.40 | 0.44 | 0.39 | 0.43 | 0.45 | 0.48 | 0.30 | 0.34 | 0.38 | 0.42 | 0.18 | 1.54 | **0.27** |
| psse0 | 0.45 | 0.64 | 0.43 | 0.63 | 0.44 | 0.64 | 0.41 | 0.62 | 0.41 | 0.62 | 0.28 | 1.01 | **0.53** |
| co9 | 0.43 | 0.47 | 0.38 | 0.42 | 0.46 | 0.50 | 0.34 | 0.39 | 0.41 | 0.46 | 0.18 | 1.54 | **0.27** |
| baxter | 0.69 | 0.75 | 0.67 | 0.74 | 0.47 | 0.57 | 0.45 | 0.56 | 0.43 | 0.54 | 0.30 | 1.09 | **0.45** |
| graphics | 0.74 | 0.87 | 0.80 | 0.91 | 0.68 | 0.84 | 0.48 | **0.75** | 0.49 | **0.75** | 0.55 | 1.01 | 0.79 |
| fome12 | 0.29 | 0.31 | 0.32 | 0.35 | 0.32 | 0.35 | 0.18 | 0.21 | 0.19 | 0.22 | 0.10 | 2.74 | **0.20** |
| route | 0.34 | 0.36 | 0.48 | 0.50 | 0.37 | 0.39 | 0.62 | 0.64 | 0.59 | 0.61 | 0.08 | 1.38 | **0.12** |
| fxm4_6 | 1.54 | 1.41 | 1.23 | 1.18 | 0.86 | 0.89 | 0.70 | **0.77** | 0.71 | 0.78 | 0.75 | 1.17 | 0.85 |
| | | | | | Geometric Means | | | | | | | | |
| Symmetric | 0.74 | 0.80 | 0.73 | 0.79 | 0.67 | 0.74 | 0.54 | 0.64 | 0.47 | 0.58 | 0.40 | 1.26 | 0.57 |
| Nonsymmetric | 0.74 | 0.76 | 0.66 | 0.70 | 0.63 | 0.68 | 0.43 | 0.51 | 0.40 | 0.48 | 0.34 | 1.54 | 0.50 |
| Rectangular | 0.60 | 0.64 | 0.57 | 0.62 | 0.51 | 0.57 | 0.41 | 0.49 | 0.43 | 0.51 | 0.27 | 1.36 | 0.40 |
| Overall | 0.69 | 0.73 | 0.65 | 0.70 | 0.60 | 0.66 | 0.45 | 0.54 | 0.43 | 0.52 | 0.33 | 1.38 | 0.48 |

TABLE 6.6
*Sensitivity of sHP$_{RN}$ [43] and sHP$_{CN}$ to cache-line size (cache size = part-weight threshold = 64 KB)*

| Line Size (Byte) | Nonsymmetric | | Rectangular | |
|---|---|---|---|---|
| | sHP$_{RN}$ | sHP$_{CN}$ | sHP$_{RN}$ | sHP$_{CN}$ |
| | $x$ | $x$ | $x$ | $x$ |
| 8 | 0.70 | 0.53 | 0.62 | 0.52 |
| 16 | 0.68 | 0.49 | 0.58 | 0.47 |
| 32 | 0.65 | 0.45 | 0.52 | 0.41 |
| 64 | 0.61 | 0.41 | 0.44 | 0.34 |
| 128 | 0.57 | 0.38 | 0.39 | 0.28 |
| 256 | 0.52 | 0.33 | 0.36 | 0.23 |
| 512 | 0.33 | 0.30 | 0.23 | 0.23 |

As seen in Table 6.6, the performance of sHP$_{RN}$ is considerably more sensitive to the cache-line size than that of sHP$_{CN}$. For nonsymmetric matrices, as the line size is increased from 8 bytes (1 word) to 512 bytes, the average normalized cache-miss count decreases from 0.70 to 0.33 in the sHP$_{RN}$ method, whereas it decreases from 0.53 to 0.30 in the sHP$_{CN}$ method. Similarly, for rectangular matrices, the average normalized cache-miss count decreases from 0.62 to 0.23 in the sHP$_{RN}$ method, whereas it decreases from 0.52 to 0.23 in the sHP$_{CN}$ method. As seen in Table 6.6, the performance of these two methods become very close for the largest line size of 512 bytes (64 words). This experimental finding conforms to our expectation that sHP$_{RN}$ exploits spatial locality better than sHP$_{CN}$, whereas sHP$_{CN}$ exploits temporal locality better than sHP$_{RN}$.

**6.2.2. Comparison of 1D and 2D methods.** We proceed with the relative performance comparison of the 1D- and 2D-partitioning-based methods which will be referred to as 1D methods and 2D methods, respectively, in the rest of the paper. As seen in Table 6.5, on the average, 2D methods sHP$_{eRCN}$ and mHP$_{RCN}$ perform better than the 1D methods sHP$_{RN}$ and sHP$_{CN}$. The performance gap between the 2D and 1D methods is considerably higher in reordering symmetric matrices in favour of 2D methods. This experimental finding may be attributed to the relatively restricted search space of the column-net model (as well as the row-net model) in 1D partitioning of symmetric matrices. The relative performance comparison of 2D methods shows that sHP$_{eRCN}$ and mHP$_{RCN}$ display comparable performance for symmetric matrices, whereas mHP$_{RCN}$ performs much better than sHP$_{eRCN}$ for nonsymmetric and rectangular matrices, on the average. mHP$_{RCN}$ performs 8.3% better than sHP$_{eRCN}$ in terms of cache misses due to the access of $x$-*vector* and $y$-*vector* entries, on the overall average.

As seen in Table 6.5, mHP$_{RCN}$ incurs significantly less $x$-*vector* entry misses than sHP$_{eRCN}$ on the overall average. This is expected because the multiple-SpMxV framework utilized in mHP$_{RCN}$ enables better exploitation of temporal locality in accessing $x$-*vector* entries. However the increase in the $y$-*vector* entry misses, which is introduced by the multiple-SpMxV framework, does not amortize in some of the reordering instances. As expected, mHP$_{RCN}$ performs better than sHP$_{eRCN}$ in the reordering of matrices that contain dense rows. For example, in the reordering of symmetric matrices `a5esindl`, `bloweya`, and `brainpc2`, which respectively contain dense rows with 9993, 10001, and 13799 nonzeros, mHP$_{RCN}$ performs significantly better than sHP$_{eRCN}$. Similar experimental findings can be observed in Table 6.5 for the following matrices that contain dense rows: square nonsymmetric matrices `circuit_4`, `ckt11752_dc_1`, `mult_dcop_01` and rectangular matrices `baxter`, `co9`, `cq9` and `route`. Although `shallow_water` and `psse0` do not contain dense rows (maximum number of nonzeros in a row is only 4 in both matrices), mHP$_{RCN}$ performs significantly better than sHP$_{eRCN}$ in reordering these two matrices. mHP$_{RCN}$ incurs significantly less cache misses due to the access of $x$-*vector* entries while incurring

TABLE 6.7

*Sensitivity of $sHP_{CN}$, $sHP_{eRCN}$, and $mHP_{RCN}$ to cache-line size (cache size = part-weight threshold = 64 KB)*

| Line Size (byte) | Single SpMxV | | | | Multiple SpMxVs | | |
|---|---|---|---|---|---|---|---|
| | $sHP_{CN}$ | | $sHP_{eRCN}$ | | $mHP_{RCN}$ | | |
| | $x$ | $x+y$ | $x$ | $x+y$ | $x$ | $y$ | $x+y$ |
| 8 | 0.59 | 0.69 | 0.59 | 0.69 | 0.47 | 1.09 | 0.62 |
| 16 | 0.55 | 0.64 | 0.55 | 0.64 | 0.43 | 1.15 | 0.58 |
| 32 | 0.50 | 0.59 | 0.49 | 0.59 | 0.36 | 1.23 | 0.52 |
| 64 | 0.45 | 0.54 | 0.43 | 0.52 | 0.33 | 1.38 | 0.48 |
| 128 | 0.41 | 0.49 | 0.38 | 0.46 | 0.28 | 1.50 | 0.43 |
| 256 | 0.37 | 0.44 | 0.33 | 0.40 | 0.24 | 1.60 | 0.38 |
| 512 | 0.36 | 0.42 | 0.30 | 0.36 | 0.25 | 1.83 | 0.38 |

very small number of additional cache misses due to the access of $y$-*vector* entries. The reason behind the latter finding is the very small number of shared rows among the $A^k$ matrices obtained by $mHP_{RCN}$ in splitting these two matrices. For example, in one of the splittings generated by $mHP_{RCN}$, among the 81920 rows of `shallow_water`, only 785 rows are shared and all of them are shared between only two distinct $A^k$ matrices.

**6.2.3. Experimental sensitivity analysis.** Table 6.7 shows the comparison of the sensitivities of the proposed methods $sHP_{CN}$, $sHP_{eRCN}$ and $mHP_{RCN}$ to the cache-line size. In the construction of the averages reported in this table, simulation results of every method are normalized with respect to those of the original ordering with the respective cache-line size. In terms of cache misses due to access of $x$-*vector* entries, the performance of each method compared to the original ordering increases with increasing cache-line size. However, in terms of cache misses due to access of $y$-*vector* entries, the performance of $mHP_{RCN}$ compared to the original ordering decreases with increasing cache-line size. So, with increasing cache-line size, the performance gap between $mHP_{RCN}$ and the other two methods $sHP_{CN}$ and $sHP_{eRCN}$ increases so that $sHP_{eRCN}$ performs better than $mHP_{RCN}$ for larger cache-line size of 512 bytes. This experimental finding can be attributed to the deficiency of the multiple-SpMxV framework in exploiting spatial locality in accessing $y$-*vector* entries. We believe that models and methods need to be investigated for intelligent global row ordering to overcome this deficiency of the multiple-SpMxV framework.

We introduce Table 6.8 to display the sensitivities (as overall averages) of the HP-based reordering methods to the part-weight threshold ($W_{max}$) used in terminating the RB process. The performance of each method increases with decreasing part-weight threshold until the part-weight threshold becomes equal to the cache size. For each method, the rate of performance increase begins to decrease as the part-weight threshold becomes closer to the cache size. The performance of each method remains almost the same with decreasing part-weight threshold below the cache size except $mHP_{RCN}$. The slight decrease in the performance of $mHP_{RCN}$ with decreasing part-weight threshold below the cache size can be attributed to the increase in the number of $y$ misses with increasing number of $A^k$ matrices because of the deficiency of the multiple-SpMxV framework in exploiting spatial locality in accessing $y$-*vector* entries. These experimental findings show the validity of Theorems 4.1, 4.2, and 5.1 for the effectiveness of the proposed $sHP_{CN}$, $sHP_{eRCN}$, and $mHP_{RCN}$ methods, respectively. Although the proposed HP-based methods are cache-size-aware methods, the ones that utilize the single-SpMxV framework can easily be modified to become cache oblivious methods by continuing the RB process until the parts become sufficiently small or the qualities of the bipartitions drop below a predetermined threshold.

**6.3. OSKI experiments.** For large size matrices, OSKI experiments are performed by running OSKI version 1.0.1h (compiled with `gcc`) on a machine with 2.66 GHz Intel Q8400 and 4 GB of RAM, where each core pair shares 2 MB 8-way set-associative L2 cache. The

TABLE 6.8
*Sensitivity of HP-based reordering methods to the part-weight threshold (cache size = 64 KB)*

| Part Size (KB) | 1D Partitioning | | | | 2D Partitioning | | | | |
|---|---|---|---|---|---|---|---|---|---|
| | sHP$_{RN}$ [43] | | sHP$_{CN}$ | | sHP$_{eRCN}$ | | mHP$_{RCN}$ | | |
| | $x$ | $x+y$ | $x$ | $x+y$ | $x$ | $x+y$ | $x$ | $y$ | $x+y$ |
| 512 | 0.79 | 0.81 | 0.71 | 0.75 | 0.69 | 0.73 | 0.63 | 1.08 | 0.69 |
| 256 | 0.68 | 0.72 | 0.61 | 0.67 | 0.57 | 0.63 | 0.49 | 1.15 | 0.58 |
| 126 | 0.62 | 0.68 | 0.51 | 0.59 | 0.48 | 0.56 | 0.39 | 1.28 | 0.52 |
| 64  | 0.60 | 0.66 | 0.45 | 0.54 | 0.43 | 0.52 | 0.34 | 1.42 | 0.50 |
| 32  | 0.59 | 0.66 | 0.43 | 0.52 | 0.42 | 0.51 | 0.33 | 1.53 | 0.51 |
| 16  | 0.60 | 0.66 | 0.43 | 0.52 | 0.42 | 0.51 | 0.34 | 1.57 | 0.52 |
| 8   | 0.61 | 0.67 | 0.43 | 0.52 | 0.42 | 0.51 | 0.35 | 1.61 | 0.54 |

Generalized Compressed Sparse Row (GCSR) format available in OSKI is used for all reordering methods. GCSR handles empty rows by augmenting the traditional CSR with an optional list of non-empty row indices thus enabling the multiple-SpMxV framework. For each reordering instance, an SpMxV workload contains 100 calls to `oski_MatMult()` with the same matrix after 3 calls as a warm-up.

Table 6.9 displays the performance comparison of the existing and proposed methods for large size matrices. In the table, the first column shows OSKI running times without tuning for original matrices. The second column shows the normalized OSKI running times obtained through the full tuning enforced by the `ALWAYS_TUNE_AGGRESSIVELY` parameter for original matrices. The other columns show the normalized running times obtained through the reordering methods. Each normalized value is calculated by dividing the OSKI time of the respective method by untuned OSKI running time for original matrices. As seen in the first two columns of the table, optimizations provided through the OSKI package do not improve the performance of the SpMxV operation performed on the original matrices. This experimental finding can be attributed to the irregularly sparse nature of the test matrices. We should mention that optimizations provided through the OSKI package do not improve the performance of the SpMxV operation performed on the reordered matrices.

The relative performance figures given in Table 6.9 for different reordering methods in terms of OSKI times in general conform to the relative performance discussions given in Section 6.2 based on the cache-miss simulation results. As seen in Table 6.9, on the overall average, the 2D methods sHP$_{eRCN}$ and mHP$_{RCN}$ perform better than the 1D methods sHP$_{RN}$ and sHP$_{CN}$, where mHP$_{RCN}$ (adopting the multiple-SpMxV framework) is the clear winner. Furthermore, for the relative performance comparison of the 1D methods, the proposed sHP$_{CN}$ method performs better than the existing sHP$_{RN}$ method. On the overall average, sHP$_{CN}$, sHP$_{eRCN}$ and mHP$_{RCN}$ achieve significant speedup by reducing the SpMxV times by 11%, 14% and 18%, respectively, compared to the unordered matrices; thus confirming the success of the proposed reordering methods.

Table 6.10 shows cache-miss simulation results for large size matrices and it is introduced to show how the performance comparison in terms of cache-miss simulations relates to performance comparison in terms of OSKI running times. In Table 6.10, the "$tot$" column shows the normalized total number of cache misses including compulsory cache misses due to the access of matrix nonzeros. The total numbers of cache misses are also displayed since these values actually determine the performance of the reordering methods in terms of OSKI times. Comparison of Tables 6.9 and 6.10 shows that the amount of performance improvement attained by the proposed methods in terms of OSKI times is in general considerably higher than the amount of performance improvement in terms of total number of cache misses. For example, sHP$_{eRCN}$ performs only 4% less cache misses then the unordered case, whereas it achieves 14% less OSKI running time, on the overall average.

Table 6.11 is introduced to evaluate the preprocessing overhead of the reordering methods. For each test matrix, the reordering times of all methods are normalized with respect to the OSKI time of the SpMxV operation using the unordered matrix and geometric averages of these normalized values are displayed in the "overhead" column of the table. In the table, the "amortization" column denotes the average number of SpMxV operations required to amortize the reordering overhead. Each "amortization value" is obtained by dividing the average normalized reordering overhead by the overall average OSKI time improvement taken from Table 6.9. Overhead and amortization values are not given for the sRCM method since sRCM does not improve the OSKI running time at all.

TABLE 6.9
*OSKI running times for large size test matrices (cache size = part-weight threshold = 2 MB)*

| | Actual | | Normalized w.r.t. Actual Times on Original Order | | | | | |
|---|---|---|---|---|---|---|---|---|
| | Original Order | | Existing Methods | | | Proposed Methods | | |
| | | | Single SpMxV | | | | | Mult. SpMxVs |
| | not tuned (ms) | OSKI tuned | sBFS [35] | sRCM [24] Modified | sHP$_{RN}$ [43] (1D Part.) | sHP$_{CN}$ (1D Part.) | sHP$_{eRCN}$ (2D Part.) | mHP$_{RCN}$ (2D Part.) |
| Symmetric Matrices | | | | | | | | |
| c-73 | 0.454 | 1.00 | 1.02 | 1.06 | 0.93 | 0.92 | 0.92 | **0.90** |
| c-73b | 0.456 | 1.00 | 1.01 | 1.07 | 0.93 | 0.92 | 0.91 | **0.89** |
| rgg_n_2_17_s0 | 0.503 | 0.95 | 0.92 | 1.07 | 0.89 | 0.82 | **0.76** | 0.91 |
| boyd2 | 0.726 | 1.19 | 1.00 | 1.14 | 0.95 | 0.92 | 0.89 | **0.85** |
| ins2 | 1.207 | 1.00 | 0.96 | 2.32 | 0.97 | 1.06 | 0.97 | **0.67** |
| rgg_n_2_18_s0 | 1.051 | 0.96 | 0.90 | 1.07 | 0.99 | 0.99 | **0.75** | 0.81 |
| Square Nonsymmetric Matrices | | | | | | | | |
| Raj1 | 0.629 | 1.04 | 0.88 | 0.96 | 0.86 | **0.82** | 0.83 | 0.84 |
| rajat21 | 0.953 | 1.07 | 1.01 | 1.16 | 1.00 | 0.95 | 0.96 | **0.90** |
| rajat24 | 0.963 | 1.02 | 1.04 | 1.16 | 0.99 | 0.94 | 0.96 | **0.91** |
| ASIC_320k | 1.436 | 0.99 | 1.09 | 1.44 | 0.97 | 0.92 | 0.73 | **0.64** |
| Stanford_Berkeley | 2.325 | 1.04 | 1.01 | 1.05 | 1.10 | 1.01 | **0.89** | 0.98 |
| Rectangular Matrices | | | | | | | | |
| kneser_10_4_1 | 0.694 | 1.02 | 0.70 | 0.87 | 0.81 | **0.67** | 0.89 | 0.68 |
| neos | 0.697 | 1.26 | 1.14 | 1.19 | 1.00 | **0.95** | 0.95 | 0.96 |
| wheel_601 | 1.377 | 1.27 | 0.82 | 0.75 | 0.69 | 0.69 | 0.66 | **0.52** |
| LargeRegFile | 2.643 | 1.55 | 1.19 | 1.30 | 1.04 | **0.95** | 0.95 | 0.96 |
| cont1_l | 2.939 | 1.14 | 1.04 | 1.19 | 1.05 | **0.93** | 0.93 | 0.95 |
| degme | 2.770 | 1.04 | 0.77 | 1.26 | 0.87 | 0.77 | 0.78 | **0.74** |
| Geometric Means | | | | | | | | |
| Symmetric | - | 1.01 | 0.97 | 1.23 | 0.94 | 0.94 | 0.86 | 0.84 |
| Nonsymmetric | - | 1.03 | 1.00 | 1.14 | 0.98 | 0.93 | 0.87 | 0.84 |
| Rectangular | - | 1.20 | 0.93 | 1.07 | 0.90 | 0.82 | 0.85 | 0.78 |
| Overall | - | 1.08 | 0.96 | 1.15 | 0.94 | 0.89 | 0.86 | 0.82 |

As seen in Table 6.11, top-down HP-based methods are significantly slower than the bottom-up sBFS method. The running times of two 1D methods sHP$_{RN}$ and sHP$_{CN}$ are comparable as expected. As also seen in the table, the 2D methods are considerably slower than the 1D methods as expected. In the column-net hypergraph model used in 1D method sHP$_{CN}$, the number of vertices and the number of nets are equal to the number of rows and the number of columns, respectively, and the number of pins is equal to the number of nonzeros. In the hypergraph model used in 2D methods, the number of vertices and the number of nets are equal to the number of nonzeros and the number of rows plus the number of columns, respectively, and the number of pins is equal to two times the number of nonzeros. That is, the hypergraphs used in 2D methods are considerably larger than the hypergraphs used in 1D methods. So partitioning the hypergraphs used in 2D methods takes considerably longer time than partitioning the hypergraphs used in 1D methods, and the running time difference becomes higher with increasing matrix density in favour of 1D methods. There exists a considerable difference in the running times of two 2D methods sHP$_{eRCN}$ and mHP$_{RCN}$ in favour of sHP$_{eRCN}$. This is because of the removal of the vertices connected by the cut row-nets in the enhanced row-column-net model used in sHP$_{eRCN}$.

TABLE 6.10
*Simulation results (misses) for large size test matrices (cache size = part-weight threshold = 2 MB)*

| | Existing Methods | | | | | | Proposed Methods | | | | | |
|---|---|---|---|---|---|---|---|---|---|---|---|---|
| | Single SpMxV | | | | | | | | | | M. SpMxVs | |
| | sBFS [35] | | sRCM [24] Modified | | sHP$_{RN}$ [43] (1D Part.) | | sHP$_{CN}$ (1D Part.) | | sHP$_{eRCN}$ (2D Part.) | | mHP$_{RCN}$ (2D Part.) | |
| | $x+y$ | tot | $x+y$ | tot | $x+y$ | tot | $x+y$ | tot | $x+y$ | tot | $x+y$ | tot |
| Symmetric Matrices | | | | | | | | | | | | |
| c-73 | 0.99 | 1.00 | 1.00 | 1.00 | 0.69 | 0.94 | 0.59 | 0.92 | 0.60 | 0.92 | **0.53** | 0.91 |
| c-73b | 0.99 | 1.00 | 1.00 | 1.00 | 0.69 | 0.94 | 0.59 | 0.92 | 0.60 | 0.92 | **0.53** | 0.91 |
| rgg_n_2_17_s0 | **0.97** | 1.00 | 1.02 | 1.00 | 1.10 | 1.01 | 1.14 | 1.01 | 1.02 | 1.00 | 0.98 | 1.00 |
| boyd2 | 1.02 | 1.01 | 1.09 | 1.14 | 0.83 | 0.94 | 0.72 | 0.91 | 0.68 | 0.90 | **0.55** | 0.85 |
| ins2 | 0.94 | 0.98 | 1.31 | 1.18 | 0.94 | 0.98 | 0.94 | 0.98 | 0.89 | 0.96 | **0.19** | 0.71 |
| rgg_n_2_18_s0 | **0.98** | 1.00 | 1.01 | 1.00 | 1.69 | 1.05 | 1.58 | 1.04 | 1.05 | 1.00 | 1.00 | 1.00 |
| Square Nonsymmetric Matrices | | | | | | | | | | | | |
| Raj1 | 0.98 | 0.99 | 0.96 | 0.99 | 0.73 | 0.94 | 0.62 | 0.91 | 0.66 | 0.92 | **0.58** | 0.90 |
| rajat21 | 1.36 | 1.08 | 1.37 | 1.08 | 1.15 | 1.03 | 0.88 | 0.97 | 0.98 | 1.00 | **0.75** | 0.95 |
| rajat24 | 1.53 | 1.11 | 1.46 | 1.09 | 1.09 | 1.02 | 0.81 | 0.96 | 0.96 | 0.99 | **0.67** | 0.93 |
| ASIC_320k | 1.61 | 1.17 | 1.61 | 1.16 | 0.89 | 0.97 | 0.79 | 0.94 | 0.57 | 0.88 | **0.32** | 0.82 |
| Stanford_Berkeley | 1.20 | 1.02 | 1.65 | 1.07 | 1.48 | 1.04 | 0.94 | 0.99 | 1.07 | 1.01 | **0.71** | 0.97 |
| Rectangular Matrices | | | | | | | | | | | | |
| kneser_10_4_1 | 1.33 | 1.09 | 1.50 | 1.13 | 1.18 | 1.05 | **0.85** | 0.96 | 1.14 | 1.04 | **0.85** | 0.97 |
| neos | 1.12 | 1.03 | 1.13 | 1.03 | 1.00 | 1.00 | **0.92** | 0.98 | **0.92** | 0.98 | **0.92** | 0.98 |
| wheel_601 | 1.39 | 1.10 | 1.40 | 1.10 | 1.16 | 1.04 | 1.02 | 1.00 | 1.14 | 1.03 | **0.91** | 0.98 |
| LargeRegFile | 1.99 | 1.20 | 1.89 | 1.18 | 1.56 | 1.11 | **1.00** | 1.00 | **1.00** | 1.00 | **1.00** | 1.00 |
| cont1_l | 1.25 | 1.06 | 1.27 | 1.07 | 1.01 | 1.00 | **0.75** | 0.94 | **0.75** | 0.94 | 0.76 | 0.94 |
| degme | 0.35 | 0.86 | 1.06 | 1.01 | 0.68 | 0.93 | 0.36 | 0.86 | 0.43 | 0.88 | **0.21** | 0.83 |
| Geometric Means | | | | | | | | | | | | |
| Symmetric | 0.98 | 1.00 | 1.07 | 1.05 | 0.94 | 0.98 | 0.87 | 0.96 | 0.78 | 0.95 | 0.55 | 0.89 |
| Nonsymmetric | 1.32 | 1.07 | 1.38 | 1.08 | 1.04 | 1.00 | 0.80 | 0.96 | 0.82 | 0.96 | 0.58 | 0.91 |
| Rectangular | 1.10 | 1.05 | 1.35 | 1.09 | 1.06 | 1.02 | 0.77 | 0.96 | 0.85 | 0.98 | 0.69 | 0.95 |
| Overall | 1.12 | 1.04 | 1.25 | 1.07 | 1.01 | 1.00 | 0.81 | 0.96 | 0.82 | 0.96 | 0.61 | 0.92 |

TABLE 6.11
*Average normalized reordering overhead and average number of SpMxV operations required to amortize the reordering overhead*

| | Existing Methods | | | | Proposed Methods | | | | | |
|---|---|---|---|---|---|---|---|---|---|---|
| | Single SpMxV | | | | | | | | Multiple SpMxVs | |
| | sBFS [35] | | sHP$_{RN}$ [43] (1D Part.) | | sHP$_{CN}$ (1D Part.) | | sHP$_{eRCN}$ (2D Part.) | | mHP$_{RCN}$ (2D Partioning) | |
| | Overhead | Amortization | Overhead | Amortization | Overhead | Amortization | Overhead | Amortization | Overhead | Amortization |
| Symmetric | 17 | 465 | 194 | 3135 | 190 | 1716 | 514 | 3732 | 920 | 5097 |
| Nonsymmetric | 26 | 700 | 314 | 5078 | 304 | 2740 | 664 | 4822 | 1198 | 6640 |
| Rectangular | 23 | 621 | 383 | 6197 | 254 | 2292 | 620 | 4503 | 1240 | 6870 |
| Overall | 22 | 587 | 286 | 4620 | 245 | 2209 | 596 | 4327 | 1110 | 6149 |

As seen in Table 6.11, the top-down HP methods amortize for larger number of SpMxV computations compared to the bottom-up sBFS method. For example, the use of sHP$_{CN}$ instead of sBFS amortizes after 276% more SpMxV computations on the overall average. As also seen in the table, 2D methods amortize for larger number of SpMxV computations compared to the 1D methods. For example, the use of mHP$_{RCN}$ instead of sHP$_{CN}$ amortizes after 178% more SpMxV computations.

**7. Conclusion.** Single- and multiple-SpMxV frameworks were investigated for exploiting cache locality in SpMxV computations that involve irregularly sparse matrices. For the single-SpMxV framework, two cache-size-aware top-down row/column-reordering methods based on 1D and 2D sparse matrix partitioning were proposed by utilizing the column-net and enhancing the row-column-net hypergraph models of sparse matrices. The multiple-SpMxV framework requires splitting a given matrix into a sum of multiple nonzero-disjoint matrices so that the SpMxV operation is computed as a sequence of multiple input- and output-dependent SpMxV operations. For this framework, a cache-size aware top-down matrix splitting method based on 2D matrix partitioning was proposed by utilizing the row-column-net hypergraph model of sparse matrices. The proposed hypergraph-partitioning (HP) based methods in the single-SpMxV framework primarily aim at exploiting temporal locality in accessing input-vector entries and the proposed HP-based method in the multiple-SpMxV framework primarily aims at exploiting temporal locality in accessing both input- and output-vector entries.

The performance of the proposed models and methods were evaluated on a wide range of sparse matrices. The experiments were carried out in two different settings: cache-miss simulations and actual runs by using OSKI. Experimental results showed that the proposed methods and models outperform the state-of-the-art schemes and also these results conformed to our expectation that temporal locality is more important than spatial locality (for practical line sizes) in SpMxV operations that involve irregularly sparse matrices. The two proposed methods that are based on 2D matrix partitioning were found to perform better than the proposed method based on 1D partitioning at the expense of higher reordering overhead, where the 2D method within the multiple-SpMxV framework was the clear winner.

**Acknowledgments.** This work was financially supported by the PRACE project funded in part by the EU's 7th Framework Programme under grant agreement number: RI-261557.

REFERENCES

[1] R. C. AGARWAL, F. G. GUSTAVSON, AND M. ZUBAIR, *A high performance algorithm using pre-processing for the sparse matrix-vector multiplication*, in Proceedings Supercomputing'92, Minn., MN, Nov. 1992, IEEE, pp. 32–41.

[2] I. AL-FURAIH AND S. RANKA, *Memory hierarchy management for iterative graph structures*, in IPPS/SPDP, 1998, pp. 298–302.

[3] C. AYKANAT, A. PINAR, AND Ü. V. ÇATALYÜREK, *Permuting sparse rectangular matrices into block-diagonal form*, SIAM Journal on Scientific Computing, 26 (2004), pp. 1860–1879.

[4] M. W. BERRY, B. HENDRICKSON, AND P. RAGHAVAN, *Sparse matrix reordering schemes for browsing hypertext*, Lectures in Applied Mathematics, 32 (1996), pp. 99–123.

[5] D. A. BURGESS AND M. B. GILES, *Renumbering unstructured grids to improve the performance of codes on hierarchical memory machines*, Technical report NA-95/06, Numerical Analysis Group, Oxford University Computing Laboratory, Wolfson Building, Parks Road, Oxford, OX1 3QD, May 1995.

[6] Ü. ÇATALYÜREK AND C. AYKANAT, *Hypergraph-partitioning-based decomposition for parallel sparse-matrix vector multiplication*, IEEE Trans. Parallel Dist. Systems, 10 (1999), pp. 673–693.

[7] Ü. V. ÇATALYÜREK, C. AYKANAT, AND B. UCAR, *On two-dimensional sparse matrix partitioning: Models, methods, and a recipe*, SIAM Journal on Scientific Computing, 32 (2010), pp. 656–683.

[8] Ü. V. ÇATALYÜREK AND C. AYKANAT, *Decomposing irregularly sparse matrices for parallel matrix-vector multiplications*, in Proceedings of 3rd International Symposium on Solving Irregularly Structured Problems in Parallel, Irregular'96, vol. 1117 of Lecture Notes in Computer Science, Springer-Verlag, 1996, pp. 75–86.

[9] ———, *Hypergraph-partitioning based decomposition for parallel sparse-matrix vector multiplication*, IEEE Transactions on Parallel and Distributed Systems, 10 (1999), pp. 673–693.

[10] ———, *PaToH: A Multilevel Hypergraph Partitioning Tool, Version 3.0*, Bilkent University, Department of Computer Engineering, Ankara, 06533 Turkey. PaToH is available at http://bmi.osu.edu/~umit/software.htm, 1999.

[11] ———, *A fine-grain hypergraph model for 2d decomposition of sparse matrices*, Parallel and Distributed Processing Symposium, International, 3 (2001), p. 30118b.

[12] J. M. CRUMMEY, D. WHALLEY, AND K. KENNEDY, *Improving memory hierarchy performance for irregular applications using data and computation reorderings*, in International Journal of Parallel Programming, 2001, pp. 425–433.

[13] R. DAS, D. J. MAVRIPLIS, J. SALTZ, S. GUPTA, AND R. PONNUSAMY, *The design and implementation of a parallel unstructured euler solver using software primitives*, in AIAA Journal, 1992.

[14] J. DEMMEL, J. DONGARRA, A. RUHE, AND H. VAN DER VORST, *Templates for the solution of algebraic eigenvalue problems: a practical guide*, Society for Industrial and Applied Mathematics, Philadelphia, PA, USA, 2000.

[15] C. DING AND K. KENNEDY, *Improving cache performance in dynamic applications through data and computation reorganization at run time*, SIGPLAN Not., 34 (1999), pp. 229–241.
[16] E. ELMROTH, F. GUSTAVSON, I. JONSSON, AND B. KGSTRM, *Recursive blocked algorithms and hybrid data structures for dense matrix library software*, SIAM Review, 46 (2004), pp. 3–45.
[17] J. D. FRENS AND D. S. WISE, *Auto-blocking matrix-multiplication or tracking blas3 performance from source code*, SIGPLAN Not., 32 (1997), pp. 206–216.
[18] G. HAASE, M. LIEBMANN, AND G. PLANK, *A hilbert-order multiplication scheme for unstructured sparse matrices*, Int. J. Parallel Emerg. Distrib. Syst., 22 (2007), pp. 213–220.
[19] H. HAN AND C. TSENG, *Exploiting locality for irregular scientific codes*, IEEE Trans. Parallel Distrib. Syst., 17 (2006), pp. 606–618.
[20] S. A. HAQUE AND S. HOSSAIN, *A note on the performance of sparse matrix-vector multiplication with column reordering*, Computing, Engineering and Information, International Conference on, 0 (2009), pp. 23–26.
[21] K. HELSGAUN, *An effective implementation of the lin-kernighan traveling salesman heuristic*, European Journal of Operational Research, 126 (2000), pp. 106 – 130.
[22] B. HENDRICKSON AND T. G. KOLDA, *Partitioning rectangular and structurally nonsymmetric sparse matrices for parallel processing*, SIAM Journal on Scientific Computing, 21 (2000), pp. 2048–2072.
[23] D. B. HERAS, V. B. PÉREZ, J. C. CABALEIRO, AND F. F. RIVERA, *Modeling and improving locality for the sparse-matrix-vector product on cache memories*, Future Generation Comp. Syst, 18 (2001), pp. 55–67.
[24] E. IM AND K. YELICK, *Optimizing sparse matrix vector multiplication on SMPs*, May 25 1999.
[25] G. JIN AND M. J. CRUMMEY, *Using space-filling curves for computation reordering*, in Proceedings of the Los Alamos Computer Science Institute, 2005.
[26] G. KARYPIS, V. KUMAR, R. AGGARWAL, AND S. SHEKHAR, *hMeTiS A Hypergraph Partitioning Package Version 1.0.1*, University of Minnesota, Department of Comp. Sci. and Eng., Army HPC Research Center, Minneapolis, 1998.
[27] J. KOSTER, *Parallel Templates for Numerical Linear Algebra, a High-Performance Computation Library*, master's thesis, Utrecht University, July 2002.
[28] T. LENGAUER, *Combinatorial Algorithms for Integrated Circuit Layout*, Willey–Teubner, Chichester, U.K., 1990.
[29] R. MIRCHANDANEY, J. H. SALTZ, R. M. SMITH, D. M. NICO, AND K. CROWLEY, *Principles of runtime support for parallel processors*, in ICS '88: Proceedings of the 2nd international conference on Supercomputing, New York, NY, USA, 1988, ACM, pp. 140–152.
[30] J. C. PICHEL, D. B. H., J. C. C., AND F. F. RIVERA, *Increasing data reuse of sparse algebra codes on simultaneous multithreading architectures*, Concurrency and Computation: Practice and Experience, 21 (2009), pp. 1838–1856.
[31] J. C. PICHEL, D. B. HERAS, J. C. CABALEIRO, AND F. F. RIVERA, *Performance optimization of irregular codes based on the combination of reordering and blocking techniques*, Parallel Computing, 31 (2005), pp. 858–876.
[32] A. PINAR AND M. T. HEATH, *Improving performance of sparse matrix-vector multiplication*, in Supercomputing '99: Proceedings of the 1999 ACM/IEEE conference on Supercomputing (CDROM), New York, NY, USA, 1999, ACM, p. 30.
[33] Y. SAAD, *Iterative Methods for Sparse Linear Systems, Second Edition*, Society for Industrial and Applied Mathematics, April 2003.
[34] M. M. STROUT, L. CARTER, AND J. FERRANTE, *Compile-time composition of run-time data and iteration reorderings*, SPNOTICES: ACM SIGPLAN Notices, 38 (2003).
[35] M. M. STROUT AND P. D. HOVLAND, *Metrics and models for reordering transformations*, in Proc. of the Second ACM SIGPLAN Workshop on Memory System Performance (MSP04), Washington DC., June 2004, ACM, pp. 23–34.
[36] O. TEMAM AND W. JALBY, *Characterizing the behavior of sparse algorithms on caches*, in Proceedings Supercomputing'92, Minn., MN, Nov. 1992, IEEE, pp. 578–587.
[37] A. D. TIMOTHY, *University of florida sparse matrix collection*, NA Digest, 92 (1994).
[38] S. TOLEDO, *Improving memory-system performance of sparse matrix-vector multiplication*, in IBM Journal of Research and Development, 1997.
[39] B. UCAR AND C. AYKANAT, *Partitioning sparse matrices for parallel preconditioned iterative methods*, SIAM Journal on Scientific Computing, 29 (2007), pp. 1683–1709.
[40] B. VASTENHOUW AND R. H. BISSELING, *A two-dimensional data distribution method for parallel sparse matrix-vector multiplication*, SIAM Review, 47 (2005), pp. 67–95.
[41] R. VUDUC, J. W. DEMMEL, AND K. A. YELICK, *Oski: A library of automatically tuned sparse matrix kernels*, Journal of Physics: Conference Series, 16 (2005), pp. 521+.
[42] J. WHITE AND P. SADAYAPPAN, *On improving the performance of sparse matrix-vector multiplication*, in In Proceedings of the International Conference on High-Performance Computing, IEEE Computer Society, 1997, pp. 578–587.

[43] A. N. YZELMAN AND ROB H. BISSELING, *Cache-oblivious sparse matrix–vector multiplication by using sparse matrix partitioning methods*, SIAM Journal on Scientific Computing, 31 (2009), pp. 3128–3154.